\documentclass[english]{smfart}

\usepackage{amsmath}
\usepackage{amscd}
\usepackage{amssymb}
\usepackage{amsfonts}
\usepackage{latexsym}
\usepackage{verbatim}
\usepackage{epsfig}
\usepackage[all]{xy}
\begin{document}

\newcommand\A{{\mathcal A}}
\newcommand\arrow{\leadsto}
\newcommand\C{{\mathcal C}}
\newcommand\D{{\mathcal D}}
\newcommand\diam{{\operatorname{diam}}}
\newcommand\eps{{\epsilon}}
\newcommand\erg{{\operatorname{erg}}}
\newcommand\expa{{\operatorname{exp}}}
\newcommand\feq{\approx_f}
\newcommand\Fix{{\operatorname{Fix}}}
\newcommand\fol{{\text{Fol}}}
\newcommand\fred{\succeq_f}
\newcommand\id{{\text{Id}}}
\newcommand\ieq{\approx_i}
\newcommand\ired{\succeq_i}
\newcommand\itree{\mathcal T_i}
\renewcommand\L{{\mathcal L}}
\newcommand\Lip{\operatorname{Lip}}
\newcommand\loc{{\operatorname{loc}}}
\newcommand\mlc{{\frak m}}
\newcommand\MR{{\mathbf M}}
\newcommand\mult{{\operatorname{mult}}}
\newcommand\N{{\mathbb N}}
\newcommand\oX{{\underline{X}}}
\newcommand\ov[1]{{\overline{#1}}}
\renewcommand\P{{\mathcal P}}
\newcommand\Prob{{\operatorname{Prob}}}
\newcommand\QED{$\square$}
\newcommand\QFT{QFT{} }
\newcommand\R{{\mathbb R}}
\newcommand\Sh{{\hat{\Sigma}}}
\newcommand\Shs{{\hat{\Sigma}}^*}
\newcommand\sQFT{$\ast$-QFT{} }
\newcommand\T{{\mathbf T}}
\newcommand\tFix{{\widetilde{\operatorname{Fix}}}}
\renewcommand\top{{\operatorname{top}}}
\newcommand\tP{\tilde{P}}
\newcommand\tr{{\text{Tr}}}
\newcommand\wloc{{\operatorname{wloc}}}
\newcommand\X{{\mathcal X}}
\newcommand\Z{{\mathbb Z}}

\newcommand\un[1]{{\underline{#1}}}
\newcommand\ignore[1]{}
\newcommand\new[1]{\emph{#1}}
\newcommand\NEW[1]{{\bf #1}}

\newcommand\annotation[1]{\marginpar{{\bf #1}}}

\newtheorem{theo}{Theorem}
\newtheorem{coro}{Corollary}[section]
\newtheorem{defi}[coro]{Definition}
\newtheorem{lem}[coro]{Lemma}
\newtheorem{prop}[coro]{Proposition}
\newtheorem{fact}[coro]{Fact}
\newtheorem{claim}[coro]{Claim}

\newenvironment{demo}
      {\medbreak\noindent{\sc Proof:}}
      {\hfill\QED\medbreak}

\newenvironment{demof}[1]
      {\medbreak\noindent{\sc Proof of {#1}:}}
      {\hfill\QED\medbreak}

\newenvironment{remk}
      {\medbreak\noindent{\sl Remark.}}
      {\medbreak}

\newenvironment{remks}
      {\medbreak\noindent{\sl Remarks.}}
      {\medbreak}

\newenvironment{exam}
      {\medbreak\noindent{\sl Example.}}
      {\medbreak}

\newenvironment{claim*}
      {\medbreak\noindent{\bf Claim.} \sl}
      {\medbreak}

\title{Puzzles of Quasi-Finite Type, Zeta Functions and Symbolic Dynamics for
Multi-Dimensional Maps}

\author{J\'er\^ome Buzzi}
\address{C.N.R.S. / Universit\'e Paris-Sud \\ 90405 Orsay cedex, France}
\email{jerome.buzzi@math.u-psud.fr}
\urladdr{www.jeromebuzzi.com}

\maketitle

\tableofcontents

 \pagestyle{myheadings} \markboth{\normalsize\sc J\'er\^ome
Buzzi}{\normalsize\sc Puzzles of Quasi-Finite Type}

\section{Introduction}

In what sense(s) can a dynamical system be ``complex'' and what is
the interplay between this complexity and the more classical
dynamical properties? A very large body of works has been devoted
to this basic question, especially to prove various forms of
complexity from dynamical assumptions. We are interested in
reversing this direction:

\bigbreak

\centerline{\it What are the dynamical consequences of complexity?}

\smallbreak

\centerline{\it Can complexity characterize a dynamical system?}

\bigbreak

This type of question has been studied mainly in low-complexity
settings (see, e.g., \cite{Berthe} and the references therein). We
have shown, first in a smooth setting, that a {\bf
high-complexity} assumption (which we called {\sl
entropy-expansion}) also has very thorough dynamical implications
\cite{EE}. A remarkable feature is that this condition, which
involves only so-called dimensional entropies, is enough to
analyze measures of maximum entropy and the related periodic
points. We are even able to classify such systems with respect to
all their ergodic and invariant measures of high entropy. Thus
complexity can be analyzed using only (simple) complexity
assumptions.

The proofs in \cite{EE} mix both combinatorial/entropic
arguments and geometric ones involving Lyapunov exponents, the
smoothness and the ensuing approximations by polynomials, raising
the question of separating completely both issues. In \cite{QFT}, we
achieved this separation for, e.g., subshifts of finite
type, piecewise monotonic interval maps with nonzero entropy,
and multidimensional $\beta$-transformations
giving a common proof to their common ``complexity" properties.
However the estimate required by \cite{QFT}  seems tractable
only when cylinders are connected, preventing until now the
application of these constructions to multi-dimensional, non-linear,
entropy-expanding maps.

The present paper overcomes this obstacle (see the remark after
Prop. \ref{prop-hC-geo}) by introducing a suitable type of
{\it symbolic dynamics} which we call {\bf puzzles of quasi-finite
type} --these are puzzles in the sense of Yoccoz' construction in
complex dynamics. In comparison with our work in
\cite{EE} we have to make an additional, probably
generic, assumption, but we also obtain more detailed information
about the periodic points.

The puzzle of quasi-finite type are the generalization of the
subshifts of quasi-finite type \cite{QFT} needed for
multi-dimensional, non-linear maps (see the end of section \ref{sec-simplicity}).
We generalize to these puzzles
all the results obtained for subshifts:
 \begin{itemize}
  \item existence of a finite number of ergodic probability
  measures maximizing the entropy;
  \item meromorphic extension of (suitably defined) Artin-Mazur
  zeta functions counting the periodic points;
  \item equidistribution of the periodic points;
  \item classification with respect to measure of large entropy.
 \end{itemize}
The meromorphic extension is deduced from a new, similar result
about Markov shifts (Theorem \ref{theo-zeta}) relating the radius of
meromorphy of some zeta functions of Markov shifts to their
entropy at infinity (see Definition \ref{defi-h-infty}). This is
of independent interest.

\subsection{Definitions}

We recall the notion of a puzzle due to Yoccoz
\cite{complex-puzzles} (closely related are the {\it tableaux} of
Branner and Hubbard).

\begin{defi}
A \NEW{puzzle} is $(V,i,f)$ (or just $V$), a set of pieces
$V=\sqcup_{n\geq0} V_n$ (a disjoint union of finite sets) and two
maps $i,f:V\setminus V_0\to V$ satisfying:
 \begin{itemize}
 \item $V_0$ has a single element;
  \item $i\circ f=f\circ i$;
  \item $i(V_{n+1})\subset V_n$;
  \item $f(V_{n+1})\subset V_n$.
 \end{itemize}
The \new{order} of a piece $v$ is $|v|$, the unique integer $n$
such that $v\in V_n$.
\end{defi}

The above combinatorial data defines a topological dynamics as
follows:

\begin{defi}\label{def-puz-dyn}
A puzzle $(V,i,f)$ defines the \new{dynamics} $F_V:X_V\to X_V$:
 \begin{gather*}
    X_V = \{ {v}\in V_0\times V_1\times\dots: \forall
    n\geq0\; i(v_{n+1}) = v_n\} \\
    F_V:(v_n)_{n\geq0}\longmapsto (f(v_{n+1}))_{n\geq0}.
 \end{gather*}
\end{defi}

\subsection{Some examples}
\label{sec-examples}

 For $v\in V$, we denote by $|v|$ the unique integer $n$
such that $v\in V_n$. It is the {\sl order} of $v$.

\subsubsection*{Subshifts are Puzzles}\label{sec-basic-example} Let $\sigma:\Sigma_+\to\Sigma_+$ be an
arbitrary one-sided subshift (i.e., a closed shift-invariant
subset of $\mathcal A^\N$ for some finite set $\mathcal A$, the
shift being $\sigma:(A_n)_{n\geq0}\mapsto(A_{n+1})_{n\geq0}$). We
are going to define a puzzle $V$ such that the
dynamics $F_V$ associated to $V$ is topologically conjugate to this subshift
$\Sigma_+$.

For each $n\geq0$, let $V_n$ be the set of words of length $n$
that appear in $\Sigma_+$, i.e., sequences $x_0\dots x_{n-1}$ such
that $y_{i+k}=x_i$ for $i=0,\dots,n-1$ for some $y\in X$ (by
convention, $V_0=\{\varnothing\}$ where $\varnothing$ is the empty
word). Define the two maps $i$ and $f$ by:
 \begin{gather*}
    i(A_1\dots A_n) = A_1\dots A_{n-1} \qquad \text{rightmost
    delete}\\
    f(A_1\dots A_n) = A_2\dots A_n \qquad \text{leftmost
    delete.}
 \end{gather*}

The conjugacy $h:X_V\to \Sigma_+$ is given by
 $
     h((A_1\dots A_{n})_{n\geq0}) = (A_{n+1})_{n\geq0}.
 $

 \subsubsection*{From dynamics to puzzles}\label{sec-dyn-to-puz}
Let $T:M\to M$ be a self-map. Let $\P_0=\{M\},\P_1,\P_2\dots$ be a
sequence of finite partitions of $M$ satisfying
 $$
    \P_n\preceq\P_{n+1} \text{ and }T^{-1}\P_n\preceq\P_{n+1}
 $$
(where $\P\preceq\mathcal Q$ means that $\P$ is less fine than
$\mathcal Q$: each element of $\P$ is a union of elements of
$\mathcal Q$; also $T^{-1}P:=\{T^{-1}A:A\in P\}$).

This data defines a puzzle as follows. Let $V$ be the disjoint
union of $V_n=\P_n$, $n\geq0$. Let $i(v)=w$ if $w$ is the element
of $\P_n$ containing $v\in\P_{n+1}$. Let $f(v)=w$ if $w$ is the
element of $\P_n$ containing $T(v)$ for $v\in\P_{n+1}$. The above
assumptions ensure that this is a well-defined puzzle.

Let us give several examples of this construction.

\medbreak \noindent $\bullet$ Let $\P$ be some finite partition
and let $\P_n=\P\vee T^{-1}\P\vee\dots\vee T^{-n+1}\P$. The
corresponding puzzle is topologically conjugate to the usual
symbolic dynamics, i.e., the left shift, $\sigma: (A_n)_{n\geq0}
\mapsto (A_{n+1})_{n\geq0}$ acting on:
 \begin{equation}\label{eq-usual-symbolic}
    \overline{\{A\in\P^\N:\exists x\in M\; \forall n\geq0\; T^nx\in
    A_n\}} \subset\P^\N.
 \end{equation}

\medbreak \noindent $\bullet$ Let $\P$ be some finite partition
and let $\P_n$ be the set of connected components of the elements
of $\P\vee T^{-1}\P\vee\dots\vee T^{-n+1}\P$. This is the form
used in complex dynamics \cite{complex-puzzles}.

\medbreak \noindent $\bullet$ The following is easy but important:

\begin{fact}\label{prop-anything}
Any continuous map $T$ on a Cantor set $K$ can be realized, up to
topological conjugacy, as a puzzle, that is, there exist a puzzle $V$ and a homeomorphism $\phi:X_V\to K$ with $\phi\circ F_V= T\circ\phi$.
\end{fact}

\begin{demo}
Let $\mathcal Q_n$, $n\geq1$, be a sequence of partitions of $K$
into closed-open sets with diameters going to zero. Let $\mathcal
P_{n+1}=\mathcal Q_{n+1}\vee \P_n\vee T^{-1}P_n$. It is then easy
to see that the dynamics of the puzzle thus defined is conjugate
to that of $T$.
\end{demo}

Thus, the dynamics of puzzles are even more diverse than
that of subshifts. For instance, they can have
infinite entropy or be without measures of maximum entropy. To get a tractable class we shall assume some form of ``simplicity".

\subsection{Notions of simplicity}\label{sec-simplicity}
We restrict ourselves to puzzles that are defined by ``few constraints" (and this will include subshifts of finite type as the special case of finitely many ``constraints"). The following choice of a notion of ``simplicity" turns out to allow a detailed analysis and more precisely enforces a close similarity to the classical properties of subshifts of finite type.

\subsubsection*{A notion of constraint: irreducibility}

Let the \new{$i$-tree} below $v\in V$ be the directed graph
$\itree(v)$ whose vertices are the $w\in V$ such that $$i^n(w):=\underbrace{i\circ\dots\circ i}_{n \text{factors}}(w)=v$$
for some $n\geq0$ and whose edges are  $u\to u'$ iff $u'=i(u)$.

\begin{defi}\label{def-irred}
A piece $v\in V\setminus V_0$ is \NEW{$f$-reducible} if the two
following conditions hold:
 \begin{enumerate}
  \item[(R1)] \label{i-fred} $f:\itree(v)\to\itree(f(v))$ is a graph isomorphism;
  \item[(R2)]  there is no $w\ne v$ such that $i(w)=i(v)$, $f(w)=f(v)$ and (R1)
  holds also for $w$.
 \end{enumerate}
Otherwise, $v$ is said to be \new{$f$-irreducible}.

\medbreak

\noindent {\bf Notations.} $v\fred^1 w$ means that $v$ is
$f$-reducible and $w=f(v)$. For $k>1$, $v\fred^k w$ is defined
inductively as $v\fred^1 f(v)$ and $f(v)\fred^{k-1} w$ (by
convention $v\fred^0 v$ for any $v\in V$). Finally $u\fred w$
means that $u\fred^k w$ for some $k\geq1$.
\end{defi}

\begin{remk}
Property (R1) was introduced by Yoccoz under the name of
``regularity'' .  It is equivalent to the following dynamical property
(here $[v]_V:=
\{x\in X_V:x_{|v|}=v\}$):
 $$ F_V:[v]_V \to [f(v)]_V \text{ is a bijection}.
 $$
In the setting of complex dynamics,  failure of (R1) is equivalent to containing critical points and is called criticality --see \cite{complex-dynamics}.

Condition (R2) seems new. It is often a consequence of (R1) --this
is the case, if, for instance, the restrictions $F_V|[v]_V$, $v\in
V$, are one-to-one.
\end{remk}

\begin{remk}
In the case where the puzzle is given by a subshift as in section
\ref{sec-basic-example}, condition (R2) of Definition
\ref{def-irred} is automatically satisfied (indeed, $f(v)=f(w)$
and $i(v)=i(w)$ imply $v=w$) and condition (R1) is equivalent to
the notion of a \emph{minimum left constraint} introduced for subshifts $\Sigma$ in \cite{QFT}:
it is a finite word $w_{-n}\dots w_0$ such that:
 \begin{multline*}
   \{x_0x_1\dots:x\in\Sigma \text{ s.t. } x_{-n}\dots x_0=w_{-n}\dots w_0\}
     \\
   \subsetneq \{x_0x_1\dots:x\in\Sigma \text{ s.t. } x_{-n+1}\dots x_0=w_{-n+1}\dots w_0\}.
 \end{multline*}
\end{remk}

One can understand the $f$-irreducible pieces as describing the
constraints in $X_V$.  The $f$-reducibility of some $v$ implies that the possible
$i$-extensions of $v$ are the same as those of $f(v)$. This is a
sort of ``local Markov property".   For instance, if all pieces are
reducible then $X_V=V_1^\N$.  A slightly less extreme example of
this phenomenon is the following Lemma proved in section \ref{sec-finite-irred}.

\begin{lem} \label{lem-finitely-many}
If a puzzle $V$ has only finitely many $f$-irreducible vertices,
then $X_V$ is (topologically conjugate to) a subshift of finite
type, i.e., a subshift of $\{1,2,...,d\}^\N$ for some $d\geq1$,
obtained by excluding a finite number of finite sequences
\cite{Walters-book}.
\end{lem}

More generally, one can expect puzzles with few  $f$-irreducible
vertices to be ``simple''. The definition below formalizes this
idea.

\subsubsection*{Entropy on the puzzle}

The puzzle $V$ will be equipped with the following
\new{combinatorial distance}: for $v\ne w$,
 $$
    d_V(v,w)= 2^{-n} \quad \text{if }n=\max\{0\leq k\leq \min(|v|,|w|):
       i^{|v|-k}(v)=i^{|w|-k}(w) \}.
 $$
Note that $d_V(i(v),i(w))\leq 2 d_V(v,w)$ and
$d_V(f(v),f(w))\leq 2d_V(v,w)$.

The corresponding metric on $X_V$ is
 $$
   d_V(x,y):=\sup_{n\geq0} d_V(x_n,y_n) = 2^{-n} \text{ such that
   } n=\min\{k\geq 0:x_k\ne y_k\} \text{ or }\infty.
 $$

Together with $f$, this induces a notion of Bowen balls in $V$:
for $v\in V$, $\eps>0$, $n\in\N$, the $(\eps,n)$-ball around $v$
is
 $$
 B(v,\eps,n):=\{w\in V : \forall 0\leq k<\min(n,|v|,|w|)\;
d_V(f^kw,f^kv)<\eps\}.
 $$

The covering number $r(\eps,n,S)$ is  the minimum number of
$(\eps,n)$-balls needed to cover $S\subset V$. We define the
\NEW{topological entropy of a sequence} $\mathcal S$ of subsets
$S_n\subset V_n$, $n\geq1$, as:
 $$
    h_\top(\mathcal S) = \lim_{\eps\to 0} \limsup_{n\to\infty}
    \frac1n\log r(\eps,n,S_n).
 $$

We recall first Bowen-Dinaburg formula for the topological entropy. The
\new{$(\eps,n)$-Bowen ball} at $x$ w.r.t.
$F_V$ (and a distance $d$ on $X_V$), is $B(x,\eps,n):=\{y\in
X_V:\forall k<n$ $d(F_V^kx,F_V^ky)<\eps\}$. The \new{topological
entropy} \cite{Walters-book} is
 $$
    h_\top(F_V) := \lim_{\eps\to0} h_\top(F_V,\eps)
    \text{ with }
    h_\top(F_V,\eps) = \limsup_{n\to\infty}
     \frac1n\log r(\eps,n,X_V)
 $$
where $r(\eps,n,S)$ is the minimum number of $(\eps,n)$-balls
necessary to cover $S$. We sometimes write $h_\top(V)$ instead of
$h_\top(F_V)$.

\smallbreak

Let $\mu$ be a probability measure $\mu$ on $X_V$ which is ergodic and
invariant under $F_V$. The Kolmogorov-Sinai entropy of $(F_V,\mu)$ can 
be defined as follows, according to Katok (we again refer to \cite{Walters-book}
for background):
 $$
    h(F_V,\mu) := \lim_{\eps\to0} h(F_V,\mu,\eps)
    \text{ with }
    h(F_V,\mu,\eps) = \limsup_{n\to\infty}
     \frac1n\log r(\eps,n,\mu)
 $$
where $r(\eps,n,\mu)$ is the minimum number of $(\eps,n)$-balls
whose union has $\mu$-measure at least $1/2$ (it can be proved that
$1/2$ can be replaced with any number in $(0,1)$ without affecting
$h(F_V,\mu)$).

\subsubsection*{Constraint entropy}

We now state our main condition on the complexity of puzzles

\begin{defi} \label{def-puz-sqft}
The \NEW{constraint entropy} of a puzzle $V=(V,i,f)$ is:
 $$
   h_\C(V) := h_\top((C_n)_{n\geq1})
 $$
where $C_n$ is the set of irreducible pieces of order $n$.

The puzzle $V$ is of {\bf $\ast$-quasi-finite type} (or, for
short, {\bf \sQFT}) if it satisfies:
 $$
    h_\C(V)< h_\top(V).
 $$
\end{defi}

This notion is essentially unrelated to that of subshifts of
weakly quasi-finite type defined in \cite{QFT}.

\subsubsection*{W-Local entropy}

\sQFT puzzles can still present complexity at arbitrarily small
scales. For example one can build \sQFT puzzles which are the
union of sequences of subshifts of finite type with equal or
increasing entropy so that they have either infinitely many or no
ergodic invariant probability of maximum entropy ---see section \ref{sec-nasty}.

Our second restriction prevents these phenomena.

\begin{defi}
$V$ being a \sQFT puzzle,  the
\new{W-local entropy} of $V$ is the defect in uniformity of $h(F_V,\mu)=\lim_{\eps\to0}
h(F_V,\mu,\eps)$ over \new{large entropy measures} (that is,
ergodic invariant probability measures with entropy close to the
supremum):
 $$
   h_\wloc(V):=\inf_{\eps>0} \sup_{\mu} h(F_V,\mu)-h(F_V,\mu,\eps)
 $$
where $\mu$ ranges over the ergodic invariant probability measures
on $X_V$ with entropy $>h_\C(V)$.
\end{defi}

\begin{remk}

Obviously, $h_\wloc(V)\leq h_\loc(F_V)$, the \emph{local entropy}
(introduced by Misiurewicz \cite{Misiurewicz}) under the name topological conditional
entropy) which bounds this defect in uniformity over \emph{all} measures.
In particular,  $h_\wloc(V)=0$ if $F_V$ is expansive, e.g., a subshift.
%
%(2) It might also be interesting to consider $h_\expa(T)$, the
%infimum of the numbers $h>0$ such that there exists $\eps>0$
%satisfying: for all ergodic and invariant probability measure
%$\mu$ with $h(T,\mu)>h$, $h(T,\mu)=h(T,\mu,\eps)$. It could be
%called \emph{expansiveness entropy} (not related to
%entropy-expansion). See also \cite{PV} and the references therein
%for some relations between various (robust) expansiveness and
%hyperbolic properties.
%
% Observe that if $\delta_{\operatorname{usc}}(T)$ is the defect in
%uppersemicontinuity of the entropy function over the set of
%invariant probability measures of $T$, then
% $$
%   h_\top(T)>h_\expa(T)+\delta_{\operatorname{usc}}(T)
% $$
%is, trivially, a sufficient condition for the existence of a
%measure of maximum entropy. This condition could be helpful
%notably for $C^r$ interval maps where there remains a gap between
%the known sufficient lower bound on entropy \cite{BuzziRuette} and
%that in the counter-examples \cite{SIM,RuetteMix}.
\end{remk}

\begin{defi}
A \NEW{\QFT puzzle} (or just \QFT) is a puzzle $(V,i,f)$ which
satisfies:
 $$
    h_\C(V)+h_\wloc(V) < h_\top(V).
 $$
\end{defi}

The notions of \QFT and \sQFT puzzles can be readily generalized
in the following fashion (adding new examples ---as this already
happens for subshifts, see \cite{QFT}). Observe that if $(V,i,f)$
is a puzzle, then so is: $(V^*,i^*,f^*)$ with $V^*=V$, $i^*=f$ and
$f^*=i$. $(V^*,i^*,f^*)$ is called the \NEW{dual puzzle}. The
dynamics of a puzzle and its dual are closely related.
Hence one could formally extends our
theorems by assuming that their assumptions hold either for the
puzzle or for its dual (see Sec. \ref{sec-ne-dual}).

\subsubsection*{Determinacy}  The analysis
of periodic points uses a further assumption. We
state it in terms of the projections, for $N=1,2,\dots$,
 $$
   i_N:V\to\bigcup_{k\leq N} V_k,\;
      v\mapsto i^{(|v|-N)^+}(v).
 $$
$i_N$ extends to a map $i_N:X_V\to \left(\bigcup_{k\leq N}
V_k\right)^\N$ in a natural way:
 \begin{equation}\label{def-iN-XV}
   i_N(x)=y \iff \forall k\geq0\; y_k=f^k(x_{n+k}).
 \end{equation}

\begin{defi}\label{def-determined}
A puzzle $V$ is \NEW{determined} if:
 $$
   u,v\fred^1 w \text{ and } i_1(u)=i_1(v)
    \implies u = v.
 $$
\end{defi}

\begin{remks}

(1) Many puzzles are determined, including: those defined by
subshifts and  those defined as in
Section \ref{sec-dyn-to-puz} under the extra assumptions: (i)
$T$ is one-to-one on the closure of each element of $\P$; (ii) for
each $x\in M$, $\lim_{n\to\infty}\diam(\P^n(x))=0$ where $\P^n(x)$ is the
unique element (if it exists) of $\P^n$ that contains $x$.

(2) There exist determined puzzles whose duals are not determined.
\end{remks}

\subsubsection*{QFT subshifts are determined QFT puzzles}

Let  $(\sigma,X)$ be an arbitrary QFT subshift. Let $V$ be the puzzle defined by $X$ as in section \ref{sec-examples}. As remarked above, $V$ is determined. Also, subshifts being expansive, $h(\sigma,\mu)=h(\sigma,\mu,\eps_0)$ for some $\eps_0>0$ depending only on the choice of the metric, not on $\mu$, so the W-local entropy is zero. As remarked after Definition \ref{def-irred}, the irreducible pieces of $V$ can be identified with the minimum left constraints of $X$. Hence the constraint entropies of the puzzle $V$ and of the subshift $X$ are the same.  Thus
 $$
     h_\C(V)+h_\wloc(V)=h_\C(X)<h_\top(V)=h_\top(V),
 $$
proving the claim.

\subsection{Structure Theorem}

Let us first recall the notion of entropy-conjugacy from
\cite{SIM}.

Let $T:X\to X$ be a Borel map. Let $h(T)$ be the supremum of the
entropy of all $T$-invariant probability measures. Recall that if
$X$ is compact and $T$ continuous then $h(T)$ is just the
topological entropy by the variational principle
\cite{Walters-book}. A subset $X_0\subset X$ is
\new{entropy-negligible} if it is Borel and satisfies:
 $$
   \tilde h(T,X_0):= \sup \{ h(T,\mu):\mu \text{ ergodic with }\mu(X_0)>0\} < h(T),
 $$
that is, $X_0$ is negligible in the usual sense for all large
entropy measures, i.e., invariant and ergodic probability measures
with entropy close to $h(T)$.

Two Borel maps $T:X\to X$ and $S:Y\to Y$ are
\NEW{entropy-conjugate}\footnote{The name entropy-conjugacy was
introduced by Bowen \cite{Bowen-non-cpt} for a similar
notion: topological conjugacy after discarding subsets 
having small dimension-like entropy.} if there exist entropy-negligible subsets
$X_0\subset X$ and $Y_0\subset Y$ and a Borel isomorphism
$\Psi:X\setminus X_0\to Y\setminus Y_0$ with $\Psi\circ T =
S\circ\Psi$. The \new{constant} of this isomorphism is
$\max(\tilde h(T,X_0),\tilde h(S,Y_0))$.

\medbreak

Recall also that a \NEW{Markov shift} (see \cite{GS} and also
\cite{FFY,Kitchens,Mauldin,RuetteClass,Sarig,Sarig1}, for
background) is the set $\Sigma(G)$ of all bi-infinite paths on a
countable or finite directed graph $G$ together with the
left-shift $\sigma$:
 $$
   \Sigma(G) = \{ x\in G^\Z:\forall n\in\Z\; x_n\to x_{n+1} \text{ on }G\}
   \text{ and }
   \sigma((x_n)_{n\in\Z})=(x_{n+1})_{n\in\Z}.
 $$
The Markov shifts defined by finite graphs are the classical
subshifts of finite type (of order $1$) ---see
\cite{LM} for an introduction to this rich classical theory.

$\Sigma(G)$ is \new{irreducible} if it contains a dense orbit
(equivalently $G$ is strongly connected). Any Markov shift
has a \new{spectral decomposition} as a union of countably many
irreducible Markov subshifts (up to wandering orbits). The period of a subset
$U$ of $\Sigma(G)$ is the greatest common divisor of all $k\geq1$ such
that $\sigma^k U\cap U\ne\emptyset$. The
\new{period} of $\Sigma(G)$ is the largest period of all non-empty open
subsets of $\Sigma(G)$.

A Markov shift is not compact unless it is a subshift of finite
type. Its topological entropy is therefore defined as explained
above for a general Borel system. Gurevi\v{c} \cite{Gurevic}
proved that in the irreducible case this entropy is just, for any
$(a,b)\in G^2$:
 $$
 h(G):=h(\Sigma(G))=\limsup_{n\to\infty} \frac1n\log\#\{v\in G^n:v_1=a,v_1\to
v_2\to\dots\ v_n=b\}.
 $$

An irreducible Markov shift is said to be \NEW{SPR} (for
stably positively recurrent \cite{GS}--- also called strongly
positively recurrent \cite{Sarig-pt}) if it admits an
entropy-maximizing probability measure $\mu$ which is
exponentially filling, i.e., for any non-empty open subset $U$ of
$X$,
 $$
   \lim_{n\to\infty} \frac1n\log \mu\left(X\setminus \bigcup_{k=0}^n \sigma^{-n}U \right)
     < 0.
 $$
Such Markov shifts are closest to being of finite type by a number
of results (see, e.g., Gurevi\v{c} \cite{Gurevic-SPR}, Sarig
\cite{Sarig}, Gurevi\v{c}-Savchenko \cite{GS} among others).

\medbreak

In Sec. \ref{sec-complete-D}, we shall associate to any puzzle a Markov shift $\Sigma(\mathcal D)$
defined by the adaptation to puzzles of the
``complete" Hofbauer diagram developed in \cite{SIM} for
subshifts.\footnote{This ``complete" variant essentially removes
``accidental" identifications, i.e., of the type $T(A)=T(B)$ where
$A$ and $B$ are distinct elements of the partition whereas
$T(A)=T(B)$ does not belong to that partition. This variant is
necessary for the precise counting of periodic orbits as we
explained in \cite{QFT} (it also simplifies the proof of the
partial isomorphism, see Section \ref{sec-struct}).}
\medbreak

Finally recall that the \new{natural extension} of a map $T:X\to
X$ is the ``smallest'' extension that is invertible, i.e., it is
$\tilde T:\tilde X\to\tilde X$ with $\tilde X:=\{x\in X^\Z:\forall
n\in\Z$ $T(x_n)=x_{n+1}\}$ and $\tilde
T((x_n)_{n\in\Z})=(Tx_n)_{n\in\Z}$.

\medbreak

 We may now state our key structure theorem:

\begin{theo}[Main Result]\label{theo-struct}
Let $V$ be a puzzle. Let $\mathcal D$ be its complete Markov diagram,
defined in Sec. \ref{sec-complete-D} below.

(1) If $V$ is \sQFT, the natural
extension of the dynamics of $V$ is entropy-conjugate with
constant at most $h_\C(V)$ to the Markov shift $\Sigma(\mathcal D)$.

(2) If $V$ is \QFT, then, for every
$H>h_\C(V)+h_\wloc(V)$, the spectral decomposition of $\Sigma(\mathcal D)$ contains
only finitely many irreducible Markov shifts with entropy $\geq H$.
Moreover these Markov shifts are SPR. More precisely their
entropies at infinity (see Definition \ref{defi-h-infty} below)
are at most $h_\C(V)+h_\wloc(V)$.

(3) If $V$ is both  \QFT and determined then, for any
$\eps>0$, there are an integer $N$ and a finite part
$\D_*\subset \D$ such that the following 
property holds.

There is a
period-preserving bijection between the periodic loops on $\D$
that meet $\D_*$ and the periodic orbits of $i_N(X_V)$ after
discarding a number $p_n^0$ of the $n$-periodic orbits of
$i_N(X_V)$ satisfying:
 $$
      \limsup_{n\to\infty} \frac1n\log p_n^0
          \leq h_\C(V)+h_\wloc(V)+\eps.
 $$
\end{theo}

The proof of this theorem is presented in Sections \ref{sec-struct}
and \ref{sec-period}.

\subsection{Dynamical consequences}

\subsubsection{Maximum measures}
The Structure Theorem gives the following, using Gurevi\v{c}'s
result on maximum measures for Markov shifts. Recall that a
\emph{Bernoulli scheme} is the shift $\sigma$ acting on the set of sequences
$\{1,\dots,s\}^\Z$ ($s$ a positive integer) endowed with
the invariant and ergodic probability measure $\mu$ defined by:
 $$
   \mu(\{\alpha:\alpha_0\dots\alpha_k=a_0\dots a_k\}) =  p(a_0)\dots p(a_n)
 $$
where $(p(1),\dots,p(s))$ is a probability vector. A \emph{finite extension}
is the product $(X,\sigma)$ with a permutation on a finite set.

For the sake of brevity, a \NEW{maximum measure} will be any
ergodic, invariant probability measure with maximum entropy.

\begin{theo}[Maximum Measures]\label{theo-maxmes}
A \QFT puzzle has at least one and at most finitely many maximum measure.

More precisely, those are in bijection with the SPR Markov subshifts
with maximum entropy and the natural extensions of these measures are
measure-preservingly isomorphic to finite extensions of Bernoulli schemes,

Moreover, the periods of the (cyclic permutations of the) measures and
those of the irreducible subshifts coincide.
\end{theo}

This follows from the Structure Theorem and
Gurevi\v{c} results for Markov shifts, as explained in Section
\ref{sec-max-meas}.

\begin{remk}
The proof that the \QFT condition implies the existence of a
maximum measure is closely related to a joint work
\cite{BuzziRuette} with S. Ruette.
\end{remk}

\subsubsection{Zeta functions}
We turn to the numbers of periodic points.

\begin{theo}[Zeta Functions]\label{theo-qft-zeta}
Assume that $V$ is a \QFT puzzle which is also determined. Fix a
large integer $N$ and consider the reduced zeta function:
 $$
    \zeta_N(z) := \exp \sum_{n\geq1} \frac{z^n}n \#\{x\in i_N(X_V): \sigma^n(x)=x\}.
 $$
$\zeta_N$ is holomorphic on $|z|<e^{-h_\top(V)}$ and has a
meromorphic extension to $|z|<e^{-h_\C(V)-h_\wloc(V)}$. Its
singularities near the circle
  $|z|=e^{-h_\top(V)}$ are exactly poles at
   $$
       e^{2i\pi k/p_i} e^{-h_\top(V)} \qquad i=1,\dots,r\quad
       k=0,\dots,p_i-1
   $$
(with multiplicities equal to repetitions in this list) where
$p_1,\dots,p_r$ are the periods of the distinct maximum measures
$\mu_1,\dots,\mu_r$.

Moreover, for each $\eps>0$, the poles of $\zeta_N(z)$ in
$|z|<e^{-h_\C(V)-h_\wloc(V)-\eps}$ are independent of $N$: for
$N',N>N(V,\eps)$, $\zeta_{N'}(z)/\zeta_N(z)$ extends to a
holomorphic function on this disk.
\end{theo}

This is, technically, the most delicate result as we have to go
from entropy estimates (which confuses very close points) to
counting (this is of course why the determinacy assumption is
required) ---see Section \ref{sec-qft-zeta}.

\begin{remk}

1. In contrast to \cite{QFT}, the lower-bound on the meromorphy
radius will be obtained using a new result about general Markov
shifts.

2. Counting the {\em projections} at level $N$ of periodic points
instead of the periodic points themselves is necessary as it not
even true that $\#\{x\in X_V:\sigma^n(x)=x\}<\infty$ for any
determined \QFT puzzle ---see Section \ref{sec:bad-zeta}.
\end{remk}

\subsubsection{Semi-local zeta functions for SPR Markov shifts}

The proof of Theorem \ref{theo-qft-zeta} relies on a similar (and new) result for
SPR Markov shifts. First, define the ``entropy at infinity":

\begin{defi}\label{defi-h-infty}
Let $G$ be a countable, oriented, irreducible graph. The
\NEW{entropy at infinity} of $G$ is:
  \begin{equation}\label{eq-h-infty}
    h_\infty(G)= \inf_{F\subset\subset G} \inf_{\mu_0>0}
     \sup\left\{ h(\sigma,\mu) : \mu([F])<\mu_0 \right\}
  \end{equation}
where $F$ ranges over the finite subgraphs of $G$ and
$[F]:=\{x\in\Sigma(G):x_0\in F\}$.
\end{defi}

\begin{remks}
(1) $h_\infty(G)=-\infty$ if $G$ is finite.

(2) $H\subset G$ implies that $h_\infty(H)\leq h_\infty(G)$ as
both are infimum over $\mu_0>0$ and $F\subset\subset G$ of $\sup\{
h(\sigma,\mu) : \mu([F])<\mu_0\}$ and $\sup \{ h(\sigma,\mu) :
\mu([F])<\mu_0$ and $\mu([G\setminus H])=0\}$, respectively.

(3) This definition was motivated by the observation of Ruette
\cite{RuetteClass} that the combinatorial quantities
considered by Gurevi\v{c} and Zargaryan \cite{GZ} were related to entropy at infinity.
In particular, $h_\infty(G)<h(G)$ iff $G$ is
SPR (see Proposition \ref{prop-ruette}).
\end{remks}

\begin{theo}\label{theo-zeta}
Let $\Sigma(G)$ be an irreducible Markov shift with finite Gurevic
entropy $h(G)$. For any finite subset $F\subset\subset G$, the
\NEW{semi-local zeta function} of $G$ at $F$:
 $$
    \zeta_F^G(z) := \exp \sum_{n\geq1} \frac{z^n}{n}
    \#\{x\in\Sigma(G):\sigma^n(x)=x \text{ and }
    \{x_0,\dots,x_{n-1}\}\cap F\ne\emptyset \}
 $$
is holomorphic on $|z|<e^{-h(G)}$ and has a meromorphic extension
to  $|z|<e^{-h_\infty(G)}$.

Moreover, for every $\eps>0$, there exists $F_0\subset\subset G$
such that, if $F,F'$ are two finite subsets with $F_0\subset
F,F'\subset\subset G$, then
 $$
   \frac{\zeta^G_{F'}(z)}{\zeta^G_{F}(z)}
    \text{ is holomorphic and non-zero on }|z|<e^{-(h_\infty(G)+\eps)}.
 $$
\end{theo}

\begin{remks}

(0) Notice that the semi-local zeta functions at a single
vertex coincide with the local zeta functions of \cite{GS} but
differ from those of \cite{Sarig1} (which have
usually a non-polar singularity at $z=e^{-h(G)}$ so has no
meromorphic extension).

(1) This result is new. In fact, even the case of where $F$ is
reduced to a single vertex had not been observed to our knowledge.

(2) The theorem is trivial if $h_\infty(G)=h(G)$, that is, if $G$
is not SPR (see Proposition \ref{prop-ruette}). In the opposite
extreme, for subshifts of finite type, i.e., $G$ finite, this
asserts that $\zeta^G_F$ extends meromorphically over
$\mathbb C$. Of course, in this case
$\zeta^G_F=\zeta^G/\zeta^{G\setminus F}$ in terms of the classical
Artin-Mazur zeta functions so the semi-local zeta function extends
in fact meromorphically over the Riemann sphere, i.e., is a
rational function.

(3) The conclusion of Theorem \ref{theo-zeta} is false for the
\emph{full zeta function} (i.e., $\zeta^G$). $\zeta^G$ is not
always defined as a formal series and, even if it is, can have
zero radius of convergence or it can have various types of
singularities (see \cite[Example 9.7]{GS}).

(4) For two finite subsets $F,H$, $\zeta^G_F(z)/\zeta^G_H(z)$ is meromorphic over
$|z|<\exp-h_\infty(G)$ but it is not necessarily holomorphic and
non-zero. If $G_n$ is the complete oriented graph on
$\{1,2,\dots,n\}$, we have $h(G_3)=\log 3$, $h_\infty(G_3)=-\infty$,
$\zeta^{G_3}(z)= \zeta^{G_3}_{G_3}(z) =1/(1-3z)$ and
$\zeta^{G_3}_{\{0\}}(z)=
\zeta^{G_3}(z)/\zeta^{G_2}(z)=(1-2z)/(1-3z)$.

(5) The maximum radius of a meromorphic extension\footnote{See 
Appendix A for formal definitions.} of the semi-local zeta functions
may be strictly larger than $\exp-h_\infty(G)$. Indeed, there are
Markov shifts for which the radius of meromorphy of the local zeta
functions varies (see Appendix A for an example where some local
zeta functions are rational and others have a finite radius of
meromorphy). One can wonder if these values and for instance their
supremum have a dynamical significance besides the obvious fact
that if $g\in G$ and $g'\in G'$ define the same local zeta
functions, the corresponding shifts are almost isomorphic in the
sense of \cite{BBG}. One would like to ``patch together" all the
partial informations provided by all the (semi) local zeta
functions.
\end{remks}

The proof of Theorem \ref{theo-zeta} relies on the generalization
of an algebraic formula decomposing the determinant of finite
matrices --see section \ref{sec-Markov-zeta}. In the special case
of a \emph{loop graph}, i.e., the disjoint union of $f_n$ loops
for each length $n\geq1$ based at a single vertex $a$ (see 
Appendix A), with $F$ reduced to $\{a\}$,
$h_\infty(G)=\limsup_{n\geq\infty} (1/n)\log f_n$ and the
determinantal formula coincides with the well-known identity
$\zeta_a(f)=(1-f_a(z))^{-1}$, where $f_a(z):=\sum_{n\geq1}
f_nz^n$, the rest of the proof following then that of \cite[Prop.
9.2]{GS}.

\subsubsection{Equidistribution of the periodic points}

The periodic points are equidistributed w.r.t. a suitable measure
of maximum entropy:

\begin{theo}[Equidistribution of periodic
points]\label{theo-equi-per} Assume that $V$ is a \QFT puzzle
which is determined. Let $\mu_1,\dots,\mu_r$ be the distinct
maximum measures and $p_1,\dots,p_r$ their periods, $p=
\operatorname{lcm} (p_1,\dots, p_r)$.

Fix a sufficiently large integer $N$ and consider, for $n\in p\Z$,
the measures:
 $$
    \mu_n^N := \sum_{x\in i_N(X_V) | \sigma^n(x)=x}
    \delta_x.
 $$
Then, in the weak star topology,
 $$
    \lim_{n\to\infty,n\in p\Z} \frac1{\mu_n^N(X_V)}\mu_n^N = \frac1{\sum_i p_i}\sum_{i=1}^r
    p_i\mu_i.
 $$
\end{theo}

This will also be a consequence of a result of Gurevi\v{c} and Savchenko
\cite{GS} for SPR Markov shifts.

\subsection{Classification of \QFT puzzles}
In the same way as QFT subshifts \cite{QFT}, \QFT puzzles can be
classified up to entropy-conjugacy by their entropy and periods.
Using the classification result \cite{BBG} obtained with Boyle and
Gomez for SPR Markov shifts, Theorem \ref{theo-struct} implies:

\begin{theo}[Classification]\label{theo-class}
The natural extension of \QFT puzzles are completely classified up
to entropy-conjugacy by the following data: the topological
entropy and the list, with multiplicities, of the periods of the
finitely many maximum measures.
\end{theo}

This gives a very precise meaning to our assertion that
complexity assumptions (defining \QFT puzzles) in fact
characterize them from the point of view of complexity.

\subsection{Smooth maps defining \QFT puzzles}\label{sec-puzzles-for-smooth}

We describe the class of smooth maps whose symbolic dynamics are
\QFT puzzles which will both provide interesting examples of such puzzles
and yield a new proof of variants of previous results \cite{EE} about the dynamics of
such maps.

\subsubsection*{Entropy-expansion}

Let $F:M\to M$ be a $C^\infty$ smooth map of a $d$-dimensional
compact manifold. The main assumption is that $F$ is {\bf
entropy-expanding}, which is defined as follows.
The {\em codimension one entropy} \cite{EE} is
 $$
     h^{d-1}(F) = \sup \{ h_\top(F,\phi([0,1]^{d-1}):\phi\in C^\infty(\R^{d-1}, M)
     \}.
 $$
Recall that $h_\top(F,\sigma)$ counts the number of
orbits starting from the not necessarily invariant set $\sigma$
--- see \cite{Walters-book}:
 $$%\begin{multline*}
   h_\top(F,\sigma) = \lim_{\eps\to 0} \limsup_{n\to\infty}
   r(\eps,n,\sigma).
 $$%\end{multline*}

 The {\em entropy-expanding condition} \cite{EE} is:
 $$
     h^{d-1}(F) < h_\top(F):=h_\top(F,M).
 $$
It is an open condition in the $C^\infty$ topology \cite{EE}.
Entropy-expanding maps form a natural class of multi-dimensional
non-uniformly expanding maps. This class includes all couplings of
interval maps, e.g., all self-maps of $[0,1]^2$ of the form:
 \begin{equation}\label{eq-coupling}
    (x,y)\mapsto (ax(1-x)+\eps y,by(1-y)+\eps x)
 \end{equation}
for $3.569...<a,b<4$ ($3.569...$ is the Feigenbaum parameter).

Indeed, $x\mapsto tx(1-x)$ maps $[0,1]$ into $[0,t/4]$ for $0\leq t\leq 4$ and has positive entropy for $t>3.569...$, so that for $\eps=0$, the above is entropy-expanding by \cite{BullSMF}. For $\eps\geq0$ small enough, the coupling (\ref{eq-coupling}) still preserves $[0,1]^2$. This coupling is finally entropy-expanding by the openness of this condition.

Such coupled interval maps are natural
examples of multi-dimensional non-uniformly expanding maps with
critical points but their ergodic theory has resisted all other
approaches up to now, despite all the results following
\cite{Viana} in the case where one of the two factors is assumed to
be uniformly expanding.

 \subsubsection*{Good partitions}

 We shall additionaly assume that there exists a {\em good partition}
$\P$ for $F$, i.e., with the following properties:

\begin{itemize}
  \item $\P$ is finite;
  \item each element of $\P$ is the closure of its interior;
  \item the boundary of each element of $\P$ is the image of a
  compact subset of $\R^{d-1}$ by a $C^\infty$ smooth map;
  \item the restriction of $f$ to the closure of any element $A$ of $\P$,
  $f|\bar{A}$, is one-to-one.
  \item for each $n\geq1$, each $\P,n$-cylinder:
   $$
     A_0\cap F^{-1}A_1\cap\dots \cap F^{-n+1}A_{n-1} \qquad
     A_i\in\P
   $$
  has only finitely many {\em almost connected components:}
  maximum subsets which cannot be split into two subsets at a
  positive distance;
  \item we have a uniform bound
 $$
    \sup_{x\in M} \#\{k\in\N:F^k(x)\in\partial\P\} < \infty.
 $$
 \end{itemize}

There are many  $C^\infty$ maps of compact manifolds which fail to have a good partition.
Indeed, it is not difficult to construct $C^\infty$ entropy-expanding maps which are
bounded-to-one on no open and dense set.

On the other hand we believe that among $C^\infty$ maps, having a good partition is generic, i.e., this property defines a subset  which contains a countable intersection of open and dense subsets. As a step in this direction, we prove the following in Appendix \ref{appendix-partition}:

\begin{prop} \label{prop-coupling}
The coupling in eq. (\ref{eq-coupling}) has a good partition for all parameters $(a,b,\eps)$ except for
at most a countable union of smooth hypersurfaces in $\R^3$.
\end{prop}

\subsubsection*{Puzzles of Good Entropy-Expanding Maps}

Given an entropy-expanding $C^\infty$ smooth map $F$ with a good partition $\P$ as above, we define the {\em
associated puzzle} to be $(V,i,f)$ with $V=\sqcup_{n\geq0} V_n$
where:
 \begin{itemize}
  \item $V_n$ is the collection of almost connected components of
$\P,n$-cylinders;
 \item $i,f:V_{n+1}\to V_n$ are the maps defined by $i(u)=v$ and
$f(u)=w$ if $v\supset u$ and $w\supset F(u)$;
 \end{itemize}
(this is a special case of the construction given in
section \ref{sec-examples}).
We shall see the:

\begin{theo}[Puzzles of entropy-expanding maps]\label{theo-hexp}
Let $T:M\to M$ be a $C^\infty$ smooth map of a $d$-dimensional
compact manifold. Assume that $T$ is entropy-expanding and admits a good
partition $\P$. Then the puzzle associated to $(T,\P)$ is \QFT and
determined. More precisely, $h_\wloc(V)=0$ and $h_\C(V)\leq
h^{d-1}(T)$.

In particular, such maps have finitely many maximum measures and,
up to the identifications given by some partition, their periodic
points define zeta functions with meromorphic extensions to
$|z|<\exp-h^{d-1}(T)$.
\end{theo}

These results are stated more precisely as Theorem \ref{theo-hexp-qft}
and its Corollaries \ref{coro-hexp-max}-\ref{coro-hexp-per}.

\begin{remk}
One can relax the assumption of smoothness to $C^r$ smoothness
with $r\geq1$ provided one strengthens the entropy-expansion
condition in the following way:
 $$
    h^{d-1}(f)+\frac{d-1}{r}\log^+\Lip(f) < h_\top(f)
 $$
and using that, according to \cite{SIM}, the left hand 
side dominates $H^{d-1}(f)$, the uniform codimension $1$ 
entropy defined in that work. We thus obtain a new existence result.
It generalizes the classical result of existence of
a maximum measure for all piecewise
monotone maps (i.e., $f:[0,1]\to[0,1]$ such that
$[0,1]=\bigcup_{i=1}^N[a_i,b_i]$ with $f|(a_i,b_i)$
is continuous and strictly monotone) with positive
topological entropy.
\end{remk}

\noindent{\bf Numbering.} All items are numbered consecutively
within each section, except for the theorems.

\medbreak

\noindent{\bf Acknowledgments.} J.-C. Yoccoz asked me about the
relationship between Hofbauer's towers and the puzzles of complex
dynamics a long time ago. M. Boyle was always encouraging at all
the stages this work went through. S. Ruette pointed out to me
Proposition \ref{prop-ruette} and the link with $h_\infty(G)$.
I am indebted to O. Sarig, especially with
regard to my results about zeta functions for Markov shifts. His
insights led me to the example given in the Appendix.

The observations of the referees have also significantly
improved the exposition of this paper (including section \ref{sec:semilocal-zeta}).

%%%%%%%%%%%%%%%%%%%%%%%%%%%%%%%%%%%%%%%%%%%%%%%%%%%%%%%%%%%%%%%%%%%%%%%%%%%%%%%%%%
\section{Further Examples }\label{sec-more-examples}

\subsection{Puzzles with finitely many irreducible vertices}\label{sec-finite-irred}
We prove Lemma \ref{lem-finitely-many}, i.e., that the dynamics of
a puzzle with finitely many irreducible vertices is (topologically
conjugate to) a subshift of finite type.

\medbreak

Let $n_0$ be the largest integer such that $V_{n_0}$ contains a
$f$-irreducible piece. Let $n$ be an arbitrary integer larger
than $n_0$. Recall the map $i_n:X_V\to V_n^\N$ from (\ref{def-iN-XV}).
To prove the lemma, it is enough to see that $i_n(X_V)$ is a subshift of finite
type (easy since $v$ reducible implies $F_V([v]_V)=[f(v)]_V=\bigcup_{w\in i^{-1}(v)}
[w]_V$) and that all the subshifts obtained for large $n$ are
topologically conjugate by the maps induced by the restrictions
$i:V_{m}\to V_n$, $m\geq n$.

Let $n>n_0$. Consider the finite graph $\Gamma_n$ whose vertices
are the elements of $V_n$ and whose arrows are defined by:
 $$
    u\to_{\Gamma_n}v \iff \exists w\in V_{n+1} i(w)=u \text{ and }
    f(w)=v.
 $$
Observe for future reference that, because of the definition of a
reducible vertex, $w$ above is uniquely determined by $u$ and $v$.

Let $\Sigma_n\subset V_n^\N$ be the subshift of finite type
defined by $\Gamma_n$. We claim that $i_n(X_V)=\Sigma_n$.

Observe first that $i_n(X_V)\subset \Sigma_n$. Indeed, for $x\in
X_V$ and $k\geq0$, $i((F_V^k{x})_{n+1})= (F_V^k{x})_n$
by definition of $X_V$ and $f((F_V^k{x})_{n+1})=(F_V^{k+1}{x})_n$
by definition of $F_V$.
Thus, $(F_V^{k}{x})_n \to_{\Gamma_n} (F_V^{k+1}{x})_n$,
and $i_n({x})\in\Sigma_n$.

We turn to the converse inclusion. Let $\alpha^0\in\Sigma_n$ for
some $n>n_0$. We are going to define inductively $\alpha^m\in
V_{n+m}^\N$, $m\geq1$, such that, for all $m\geq0$, $j\in\N$ and
$0\leq k\leq m$,
 \begin{equation}\label{eq-alpha}
   \text{(i) }i^k(\alpha^{m}_j)=\alpha^{m-k}_j, \text{ (ii) }
   \alpha^m_j \fred^k \alpha^{m-k}_{j+k}, \text{ (iii) }
   \alpha^m_{j+1}\in f(i^{-1}(\alpha^m_j)).
 \end{equation}
This will imply that $i_n(X_V)\supset\Sigma_n$.
Indeed, recall that $x=i_n(\alpha^0)$ means that
$x_k=i^{n-k}(\alpha_0^0)$ for $k\leq n$ and
$x_k=\alpha^{k-n}_0$ for $k\geq n$. Hence $x\in X_V$ by (i) and
$i_n(x)=\alpha^0$, as $f^k(\alpha^k_0)=\alpha^0_k$ by (ii).

Observe that (\ref{eq-alpha}) holds for $m=0$ because of the
definition of $\Gamma_n$. Let us assume that $\alpha^p_j$ has been
defined for $p\leq m$ and all $j\in\N$ so that eq.
(\ref{eq-alpha}) is satisfied. For $j\in\N$, let us build
$\alpha^{m+1}_j$ satisfying (\ref{eq-alpha}).

\begin{figure}
$$
 \xymatrix{
 \vdots & \vdots & \vdots\\
 \alpha^m_j \ar[ur]^f \ar[u]^i & \alpha^m_{j+1} \ar[u]^i \ar[ur]^f & \alpha^m_{j+2} \ar[u]^i \\
 \alpha^{m+1}_j \ar[ur]^f \ar[u]^i & \alpha^{m+1}_{j+1} \ar[u]^i \ar[ur]^f & \alpha^{m+1}_{j+2} \ar[u]^i \\
 \vdots \ar[ur]^f \ar[u]^i & \vdots \ar[u]^i \ar[ur]^f  & \vdots \ar[u]^i\\
 }
 $$
\caption{Construction of $\alpha^{m+1}$.}\label{fig-alpha}
\end{figure}
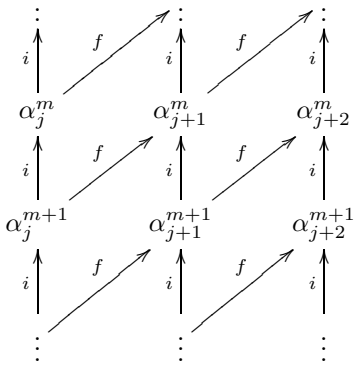

Let $\alpha^{m+1}_j\in f^{-1}(\alpha^m_{j+1})\cap
i^{-1}(\alpha^m_j)$ (this intersection is not empty by (iii), eq.
(\ref{eq-alpha}) and it is unique because $\alpha^{m+1}_j$ is
$f$-reducible). Let us check eq. (\ref{eq-alpha}) for $m+1,j$.
$i^{k+1}(\alpha^{m+1}_j)=i^k(\alpha^m_j)$ hence (i) is satisfied.
$f(\alpha^{m+1}_j)=\alpha^m_{j+1}$ and $\alpha^{m+1}_j$ is
$f$-reducible by the main assumption. Thus
$\alpha^{m+1}_j\fred^k\alpha^{m+1-k}_{j+k}$ for $k=1$ and for
$1<k\leq m$ by the induction hypothesis. This is (ii). As
$\alpha^{m+1}_j\fred^1\alpha^m_j$, the $i$-tree below
$\alpha^{m+1}_j$ is mapped by $f$ onto the $i$-tree below $\alpha^m_j$.
This gives (iii), completing the induction.

\medbreak

 Finally, one observes that $\alpha^{p+1}=i_{p+1}(\un{x})$ is uniquely
defined by $\alpha^p=i_{p}(\un{x})$ so that the natural projection
$i_{p+1}(X_V)\to i_p(X_V)$ is in fact a homeomorphism. This finishes
the proof of Lemma \ref{lem-finitely-many}.

\subsection{\sQFT puzzles with nasty dynamics}\label{sec-nasty}
We give examples of \sQFT puzzles with infinitely or no maximum
measures.

Let $\Sigma_0=\{0^\infty\},\Sigma_1,\Sigma_2,\dots$ be a sequence
of subshifts of finite type over disjoint alphabets. Assume that
the Markov order of $\Sigma_n$ is at most $n$ (i.e.,
$A\in\Sigma_n$ iff $A_k\dots A_{k+n-1}$ is a word in $\Sigma_n$
for all $k\geq0$) and that $h_\top(\Sigma_n)>0$ for all $n\geq1$.
We are going to build a puzzle which is conjugate to
$\bigcup_{n\geq0} \Sigma_n$. Taking $h_\top(\Sigma_n)=\log 2$ for
all $n\geq1$, or $h_\top(\Sigma_n)\nearrow\log 2$ as $n\to\infty$,
shall yield the required examples.

Let $L_n(\Sigma_k)$ be the set of words of length $n$ appearing in
$\Sigma_k$. The puzzle will be $(V,i,f)$ defined as follows:

Let $V_0=\{\emptyset\}$ (the empty word) and $V_n=\sqcup_{0\leq
k\leq n} L_n(\Sigma_k)$.

Let $w:=A_1\dots A_n\in V_n$. If $w\in L_n(\Sigma_n)$, then
$f(w)=i(w)=0^{n-1}$. Otherwise, let $f(A_1\dots A_n)=A_2\dots A_n$
and $i(A_1\dots A_n)=A_1\dots A_{n-1}$.

The only vertices of $V_n$ that can  be irreducible are
those $w\in L_n(\Sigma_n)$ which are mapped by $f$ to $0^{n-1}$.
For $n>N$, all these vertices are confused with $0^n$ by $i_N$. Thus at a given level
$N$, the number of distinguishable irreducible vertices in $V_n$ is
bounded independently of $n$ so that $h_\C(V)=0$. Thus $V$ is indeed
a \sQFT puzzle.

\subsection{\QFT puzzles with bad zeta functions}\label{sec:bad-zeta} Let us describe a
 determined \QFT puzzle with infinitely many periodic orbits of any given
length so that the zeta function defined from the periodic points
(and not their projections) is not even well-defined as a formal
series.

Pick a sequence of positive integers $p_1,p_2,\dots$ such that
$\#p^{-1}(k)=\infty$ for all $k\geq1$. Modify the previous
construction taking $\Sigma_0:=\{0,1\}^2$ and, for $n\geq1$, $\Sigma_n:=\{\sigma^j\omega^n:0\leq j<p_n\}$, a periodic orbit of length $p_n\geq 1$. Take the symbols $(\omega^n_j)_{n\geq 1, 0\leq j<p_n}$ pairwise distinct and disjoint from $\{0,1\}$.

Observe that, for any word $w$ from some $\Sigma_n$, $n\geq1$, $\itree(w)$ is a linear graph whereas $\itree(w)$ for all the other words are not linear. It follows that the irreducible pieces are: (i) the one-letter words $0$ and $1$;
(ii) the words from $\Sigma_n$ of length $n$. Thus, $h_\C(V)=h_\wloc(V)=0$.

Let $u,u'$ be words, $u$ not from $\Sigma_0$ such
that $u,u'\fred^1w$ and $i_1(u)=i_1(u')$. $\itree(u)$ is then a linear graph forcing $\itree(u')$ to be
so. It follows that $u=\omega^n_j\dots\omega^n_{j+\ell}$ and $u'=\omega^m_k\dots\omega^m_{k+\ell}$ for some integers $j,k,\ell,n,m$ with $\ell\geq1$ by reducibility. Hence $f(u)=f(u')$ implies that $\omega^n_{j+1}=\omega^m_{k+1}$. By the choice of pairwise distinct symbols,
this yields $u=u'$: $V$ is determined.

Therefore $V$
is indeed a determined \QFT puzzle.

\begin{remk}
Obvious adaptations of this construction yield examples with
arbitrary growth rates of the number of periodic orbits.
\end{remk}

\section{Basic properties}\label{sec-basic}

\subsection{Some properties of $f$-reducibility}
\begin{lem}\label{lem-u-unique}
If $i(u)=i(u')$, $f^k(u)=f^k(u')$ and $u\fred^k v$ and
$f:\itree(f^{l-1}(u'))\to\itree(f^{l}(u'))$, $l=1,\dots,k$ are graph
isomorphisms then $u=u'$. In particular,
 \begin{equation}\label{eq-particular}
   i(u)=i(u') \text{ and }u\fred^k w \text{ and } u'\fred^l w
     \implies k=l \text{ and }u=u'.
 \end{equation}
\end{lem}

\begin{demo}
(\ref{eq-particular}) clearly follows from the first claim. For
$k=0,1$, this claim follows from the definition of $\fred^k$.
Assume the claim for some $k-1\geq 0$ and let $u,u'$ and $v$ be as
in the claim for $k$. Now, $i(f^{k-1}(u))=
f^{k-1}(i(u))= f^{k-1}(i(u'))=i(f^{k-1}(u'))$ and both
$f^{k-1}(u)\fred^1 v$ and $f:\itree(f^{k-1}(u'))\to\itree(v)$ is
an isomorphism. This implies that $f^{k-1}(u)=f^{k-1}(u')=:w$ by
the definition of $\fred^1$. Now $i(u)=i(u')$ and $u\fred^{k-1}w$
and $f:\itree(f^{l-1}(u'))\to\itree(f^{l}(u'))$, $l=1,\dots,k-1$,
so the induction hypothesis implies $u=u'$.
\end{demo}

\begin{lem}\label{lem-red-prop}
If $i(u)\fred^k i(v)$ and $f^k(u)=v$ with $|u|,|v|\geq1$ and
$k\geq0$, then $u\fred^k v$. In particular, if $u$ with $|u|>1$ is
$f$-irreducible, then so is $i(u)$.
\end{lem}

\begin{demo}
$i(u)\fred^k i(v)$ implies that the $i$-trees below $i(u)$ and
$i(v)$ are isomorphic through $f^k$. This
implies the same for the sub-$i$-trees below $u$ and $v$. Assuming
by contradiction that $u\not\fred^k v$ we obtain that there exists
$w\in i^{-1}(i(u))$, $w\not= u$ with $f^k(w)=f^k(u)$, but this
would contradict that $f^k|\itree(i(u))$ is one-to-one.
\end{demo}

\subsection{Natural extension and duality}\label{sec-ne-dual}

Except in trivial cases, the dynamics $F_V:X_V\to X_V$ is
non-invertible. To obtain an invertible dynamical system, one goes
to the {\em natural extension}. It can be described as
$(\oX_V,F_V)$ with:
 \begin{gather*}
   \oX_V=\{ (v_{n,p})_{n,p}: \forall(n,p)\in\N\times\Z\;
       i(v_{n+1,p})=f(v_{n+1,p-1})=v_{n,p}\in V_n \}\\
   F_V: (v_{n,p})_{n,p} \longmapsto
       (v_{n,p+1})_{n,p}.
 \end{gather*}
The distance on $\oX_V$ is defined as:
 $
   d(x,y) = \sum_{n\geq0} 2^{-n}d_V(x_{-n},y_{-n}).
 $

Remark that $(\oX_V,F_V)$ is homeomorphic to the usual realization
of the natural extension: $\{ (\un{v}_{p})_{p\in\Z}\in X_V^\Z:
\forall p\in\Z\; F_V(\un{v}_{p-1})=\un{v}_{p} \}$.

\medbreak

The symmetry of the roles of $i$ and $f$ gives rise to a {\bf
duality} between puzzles: just exchange the maps $i$ and $f$
associated to a puzzle $(V,i,f)$. We denote by $(V^*,i^*,f^*)$ the
resulting puzzle. The natural extension of their dynamics
$F_{V}$ and $F_{V^*}$ are inverse
of each other, as the description of the natural extensions given
above makes it obvious.

\begin{remk}
As it was already the case for subshifts of quasi-finite type \cite{QFT},
$h_\C(V^*)$ may be different from $h_\C(V)$. It may indeeed occur that
$h_\C(V)<h_\top(V)$ and $h_\C(V^*)=h_\top(V^*)$ (or the other way around).
This allows easy construction of puzzles such that $h_\wloc(V^*)$ is
different from $h_\wloc(V)$.
%On the other hand, if $\mu$ is an invariant probability measure of $(X_V,F_V)$,
%$\un{\mu}$ the corresponding one for the natural extension and $\mu^*$
%the one for the dual puzzle dynamics,
%  $$
%     h(F_{V^*},\mu^*,2\eps) \leq h(F_V,\un{\mu},2\eps) \leq h(F_V,\mu,\eps) \leq  h(F_V,\un{\mu},\eps) \leq  h(F_{V^*},\mu^*,\eps/2)
%  $$
%for all $\eps>0$. To see, e.g., the second inequality above, observe that, letting $n_0:=\log \eps^{-1}/\log 2$:
% $$
%   [ \forall 0\leq k<n \; d(x_k,y_k)<\eps ] \implies
%       [ \forall n_0\leq k<n-n_0 \; d(F_V^kx,F_V^ky)<2\eps ] .
% $$
\end{remk}

%%%%%%%%%%%%%%%%%%%%%%%%%%%%%%%%%%%%%%%%%%%%%%%%%%%%%%%%%%%%%%%%%%%%%%%%%%%%%%%%%%
\section{Measure-theoretic Structure}\label{sec-struct}
%%%%%%%%%%%%%%%%%%%%%%%%%%%%%%%%%%%%%%%%%%%%%%%%%%%%%%%%%%%%%%%%%%%%%%%%%%%%%%%%%%

In this section we begin the proof of the structure theorem (Theorem \ref{theo-struct}). We first
introduce the Markov shift which underlies our analysis and then we
explain its consequences for entropy-conjugacy. The
proof then has three stages: (i) the Markov shift is shown to be
measurably conjugate to a part of the natural extension of the
puzzle dynamics; (ii) the entropies of the measures living on the
excluded part are bounded, yielding claim (1) of the Theorem; (iii) the entropy ``at infinity'' in
the Markov diagram is also controlled, yielding claim (2) of the Theorem. Claim (3), on the periodic points, is proved
in the next section.

\subsection{The complete Markov diagram}\label{sec-complete-D}
The key object is the following countable oriented graph.

\begin{defi}
 The {\bf complete Markov diagram} of a puzzle $V$ is a
countable, oriented graph $\mathcal D$ defined as follows. Its
vertices are the $f$-irreducible vertices of $V$. Its arrows are
the following:
 \begin{equation}\label{eq-md}
  v\arrow w \iff \exists u\in V \; i(u)= v \text{ and } u\fred w.
 \end{equation}
\end{defi}

Notice that because of Lemma \ref{lem-u-unique}, $u$ in eq.
(\ref{eq-md}) is unique given $v\arrow w$.

\begin{remk}
If $V$ is in fact a subshift over alphabet $V_1$, this complete Markov diagram reduces to the one introduced in \cite{SIM}. Under the additional assumption that  there are no
``accidental" identifications, i.e., $F_V^{|v|}([v])=F_V^{|w|}([w])$ only if $w=f^{|v|-|w|}(v)$ (assuming $|v|\geq |w|$), this further reduces to the Hofbauer diagram \cite{Hofbauer0}.
\end{remk}

\medbreak

Let $\Sigma_+(\D)$ be the associated one-sided subshift:
 $$
    \Sigma_+(\D) = \{ \un{v}\in V^\N: v_0\arrow v_1\arrow v_2\arrow\dots
    \}
 $$
together with the left-shift $\sigma((v_n)_{n\in\N})=(v_{n+1})_{n\in\N}$.

\medbreak

We build a conjugacy from the Markov shift onto (a part of) the
puzzle dynamics.

\begin{prop}\label{prop-build-wm}
Let $\un{v}\in\Sigma_+(\mathcal D)$ and $n\geq0$. There exists a
unique $w^{(n)}\in V$ such that:
 \begin{enumerate}
  \item[(i)] $i^{n}(w^{(n)})=v_0$;
  \item[(ii)] for all $k=0,\dots,n$:
 $
     i^{k}(w^{(n)}) \fred v_{n-k}.
 $
 \end{enumerate}
Moreover, the following property holds:
 \begin{enumerate}
  \item[(iii)] $i(w^{(n+1)})=w^{(n)}$.
 \end{enumerate}
\end{prop}

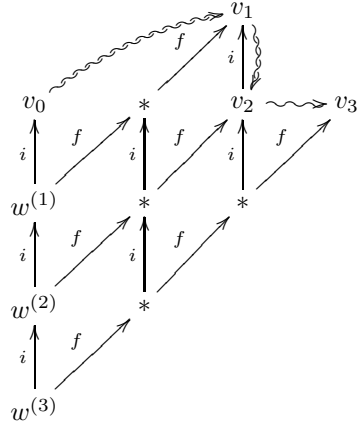
\begin{figure}
$$
 \xymatrix{
   &  &  v_1  \ar@{~>}@/^/[d] & \\
    v_0 \ar@{~>}@/^/[urr] & \ast  \ar[ur]^f &   v_2 \ar[u]^i \ar@{~>}[r] &  v_3 \\
 w^{(1)} \ar[ur]^f \ar[u]^i & \ast \ar[u]^i \ar[ur]^f & \ast \ar[u]^i \ar[ur]^f \\
 w^{(2)} \ar[ur]^f \ar[u]^i & \ast \ar[u]^i \ar[ur]^f &  \\
 w^{(3)} \ar[ur]^f \ar[u]^i &  & \\
 }
 $$
\caption{Construction of $w^{(3)}$ from $v_0v_1v_2v_3$ as in
Proposition \ref{prop-build-wm}.}\label{fig-wn}
\end{figure}

Figure \ref{fig-wn} gives a typical example of the construction of
$w^n$.

\begin{demo}
Let $\un{v}\in \Sigma_+(\mathcal D)$. For each $n\geq0$, we are
going to define $w^0,\dots,w^n$ such that:
 \begin{equation}\label{eq-mlc}
   i^j(w^j)=v_{n-j} \text{ and } \forall k=0,\dots,j\quad
       i^{k}(w^{j}) \fred v_{n-k}.
 \end{equation}
Observe that $w^{(n)}:=w^n$ will then have the required properties
(i) and (ii) by eq. (\ref{eq-mlc}). (iii) will follow from showing the
uniqueness of the solution to (\ref{eq-mlc}).

For $j=0$, set $w^{0}=v_{n}$. For $1\leq j\leq n$, assume that
$w^{j-1}$ has been defined satisfying (\ref{eq-mlc}). As
$v_{n-j}\arrow v_{n-j+1}$, there exist an integer $l\geq1$ and
$u\in i^{-1}(v_{n-j})$ such that $u\fred^l v_{n-j+1}$ (where,
necessarily, $l=|u|-|v_{n-j+1}|= |v_{n-j}|+1-|v_{n-j+1}|$). Hence
there exists a $w^{j}\in\itree(u)$ which is the $f^l$-preimage of
$w^{j-1}$ in $\itree(v_{n-j+1})$ (recall that
$i^{j-1}(w^{j-1})=v_{n-j+1}$).

Let us check (\ref{eq-mlc}) for $w^j$. Compute
 $$
   |w^j|=l+|w^{j-1}| = |v_{n-j}|+1-|v_{n-j+1}|+(j-1+|v_{n-j+1}|)
    = |v_{n-j}|+j.
 $$
As $w^j\in\itree(u)$ and $i(u)=v_{n-j}$ (i.e.,
$w^j\in\itree(v_{n-j})$), this implies the first part of
(\ref{eq-mlc}):
 \begin{equation}\label{eq-mlc1}
    i^j(w^j)=v_{n-j}.
 \end{equation}

For the second part, observe that $u\fred^l v_{n-j+1}$,
$i^{j-1}(w^j)=u$, $i^{j-1}(w^{j-1})=v_{n-j+1}$  and
$f^l(w^j)=w^{j-1}$. Hence Lemma \ref{lem-red-prop} implies that,
for $0\leq k\leq j-1$, $i^k(w^j)\fred^l i^k(w^{j-1})$. Using the
second part of (\ref{eq-mlc}) for $w^{j-1}$ we see that:
  $$
    \forall 0\leq k<j \quad i^k(w^j)\fred i^k(w^{j-1})\fred v_{n-k}.
  $$
Thus eq. (\ref{eq-mlc}) holds for $w^j$ and $k<j$. For $k=j$, this
second part is just (\ref{eq-mlc1}).

This completes the inductive construction of $w^n$.

\medbreak

For future reference, observe that $w^j$ depends only on
$v_{n-j}\dots v_n$ and that the case $k=0$ of the previous
equation gives:
 \begin{equation}\label{eq-wjfred}
    w^j\fred w^{j-1}\fred v_n
 \end{equation}

 \medbreak

 Let us check the uniqueness of problem (\ref{eq-mlc}).
We prove that for $w^n$ satisfying eq. (\ref{eq-mlc}),
$i^{n-p}w^n$ is unique by an induction on $0\leq p\leq n$. For
$p=0$, this is obvious. Assume it for $p-1\geq0$. Observe that
$i(i^{n-p}(w^n))=i^{n-p+1}(w^n)$ and $i^{n-p}(w^n)\fred v_{p}$.
These two conditions uniquely determine $i^{n-p}(w^n)$ according
to Lemma \ref{lem-u-unique}. Thus $w^{n}$ is indeed unique.

\medbreak

Thus we have shown the existence of $w^{(n)}$ satisfying
properties (i)-(iii) of the statement. We show that $w^{(n)}$ is
unique under (i) and (ii). We proceed by induction on $n$.
For $n=0$ this is obvious. Assume the uniqueness for $n-1\geq0$.
Let $w':=i(w^{(n)})$. Observe that
 \begin{itemize}
   \item $i^{n-1}(w')=i^n(w^{(n)})=v_0$;
   \item for $0\leq k<n$, $i^k(w')=i^{k+1}(w^{(n)})\fred
   v_{n-k-1}$.
 \end{itemize}
By the induction hypothesis, $w'=w^{(n-1)}$. Thus
$i(w^{(n)})=w^{(n-1)}$ and $w^{(n)}\fred v_n$. Lemma
\ref{lem-u-unique} gives the uniqueness of $w^{(n)}$,
completing the induction.
%\medbreak
% For the claim $i(w^{(n+1)})=w^{(n)}$, observe that
%$\omega:=i(w^{(n+1)})$ satisfies $i^n(\omega)=v_0$ and $i^k\omega
%\fred v_{n+1-k-1}$ for $0\leq k \leq n$. Hence, the above
%uniqueness does imply that $\omega=w^{(n)}$, i.e.,
%$i(w^{(n+1)})=w^{(n)}$ as claimed.
\end{demo}

\begin{coro}\label{coro-proj}
Let $\un{v}\in\Sigma_+(\D)$. Then there exists a unique ${x}\in X_V$
such that $x_{|v_0|}=v_0$ and for all $j\geq0$, $x_{|v_0|+j}\fred
v_j$. Moreover $x_{|v_0|+j}$ depends only on $v_0v_1\dots v_j$.
For future reference we denote this $x\in X_V$ by $x(\un{v})$.
\end{coro}

\begin{demo}
For each $n\geq0$, apply the above proposition to the sequence
$v_0\arrow\dots\arrow v_n$ to get $w^{(n)}$. As $i(w^{(n+1)})=w^{(n)}$,
we define a sequence $x$ in $X_V$ by $x_{|w^{(n)}|}=w^{(n)}$. Moreover,
for each $n\geq0$, $x_{|w^{(n)}|-n+j}=i^{n-j}(w^{(n)})\fred v_j$. As
$|w^{(n)}|=|v_0|+n$, this implies that $x_{|v_0|+j}\fred v_j$.

The uniqueness is proved by applying inductively Lemma
\ref{lem-u-unique} to $i(x_{|v_0|+j+1})=x_{|v_0|+j}$ and
$x_{|v_0|+j}\fred v_j$.

 It is obvious that $x_{|v_0|+j}=x_{|w^j|}$ depends only on
$v_0\dots v_j$.
\end{demo}

Let us define $\pi:\Sigma_+(\D)\to X_V$ by:
 $$
   \pi(\un{v}) = F_V^{|v_0|}(x(\un{v}))
 $$
with $x(\un{v})$ defined as in the above Corollary.

\begin{lem}\label{lem-D-to-X}
The map $\pi:\Sigma_+(\D)\to X_V$ is well-defined, continuous and
satisfies: $\pi\circ\sigma=F_V\circ\pi$.
\end{lem}

\begin{demo}
The above Corollary shows that $\pi$ is indeed well-defined and continuous
with values in $X_V$. We turn to the commutation relation. We must show that
 \begin{equation}\label{eq:claim-conj}
    f((\pi v)_{n})=(\pi \sigma v)_{n-1}
 \end{equation}
for all large $n$.

Let $w^{(n+1)}$ be built as in Proposition \ref{prop-build-wm}
from $v_0\dots v_{n+1}$ using the finite sequence $w^0,\dots,w^{n+1}$
defined in (\ref{eq-mlc}). Let $\tilde w^{(n)}$ be defined similarly
from $v_1\dots v_{n+1}$ using $\tilde w^0,\dots,\tilde w^n$.

Observe that $w^k=\tilde w^k$ for $k\leq n$ as they are both determined
by $v_{n+1-k}\dots v_{n+1}$. According to (\ref{eq-wjfred}),
$w^{(n+1)}\fred^\ell \tilde w^{(n)}$ with $\ell:=|w^{(n+1)}|-|\tilde w^{(n)}|
=|v_0|+1-|v_1|\geq0$. Hence
 $$
   f(f^{|v_0|}(w^{(n+1)})) = f^{|v_1|}(\tilde w^{(n)}).
 $$

Now $\pi(v)=F_V^{|v_0|}(x)$ with $x_{|v_0|+n+1}=w^{(n+1)}$. Likewise,
$\pi(\sigma(v))=F_V^{|v_1|}(y)$ with $y_{|v_1|+n}=\tilde w^{(n)}$.
(\ref{eq:claim-conj}) now follows from the previous equation.
\end{demo}

\subsection{Partial conjugacy}

We are going to show that $\pi$ gives an isomorphism between a
subset of the natural extension $\oX_V$ of $X_V$ and the whole of
$\Sigma(\D)$.

Observe that $\pi:\Sigma_+(\D)\to X_V$ extends naturally to
$\pi:\Sigma(\D)\to\oX_V$ by setting $\pi(v)=x$ with
$x_{0p}x_{1p}\dots= \pi(v_pv_{p+1}\dots)$ because of the
commutation in Lemma \ref{lem-D-to-X}.

\begin{defi}
$x\in\oX_V$ is \NEW{eventually Markovian} at time $p$ if there
exists $0\leq N<\infty$ such that:
 $$
   \forall n\geq N\quad x_{n,p-n} \fred x_{N,p-N}.
 $$
The \NEW{eventually Markovian subset} $\oX_V^M$ of $\oX_V$ is
 $$
   \oX_V^M = \{x\in\oX_V: x
    \text{ is eventually Markovian at all times }\}.
 $$
\end{defi}

\begin{prop} \label{prop-partial-iso}
Define $\iota:\oX_V^M\to\Sigma(\D)$ by $\iota(x)=v$ if, for all
$p\in\Z$, $v_p$ is the unique irreducible vertex such that for all
sufficiently large $n$:
 \begin{equation}\label{eq-iso}
x_{n,p-n}\fred v_p.
 \end{equation}
Then $\iota:(\oX_V^M,F_V)\to(\Sigma(\D),\sigma)$ is well-defined
and gives an isomorphism whose inverse is $\pi$.
\end{prop}

\begin{demo}
Let us first check that $\iota$ is well-defined with
$\iota(\oX_V^M)\subset\Sigma(\D)$. Let $x\in\oX_V^M$. $\oX_V^M$
is precisely defined so that $v=\iota(x)$ is a well-defined
element of $\D^\Z$. (\ref{eq-iso}) gives uniqueness at once.

Let us show that $v_p\leadsto v_{p+1}$ for an arbitrary $p\in\Z$.
For $n$ large enough,
 \begin{equation}\label{eq-iota}
   \text{(i) }x_{n,p-n} \fred^\ell v_p \text{ and }
   \text{(ii) }x_{n+1,p+1-n-1}=x_{n+1,p-n}\fred^k v_{p+1}
 \end{equation}
where $\ell=|x_{n,p-n}|-|v_p|$ and $k=|x_{n+1,p-n}|-|v_{p+1}|$.

Let $u=f^\ell(x_{n+1,p-n})$ (note that $|x_{n+1,p-n}|=n+1>n>\ell$).
We have $i(u)=f^\ell(i(x_{n+1,p-n}))=f^\ell(x_{n,p-n})=v_p$. Hence it is
enough to see that $u\fred v_{p+1}$. Given (\ref{eq-iota},ii),
this will follow from $\ell\leq k$. If $\ell>k$,
$x_{n+1,p-n}\fred^\ell u$ (a consequence of $x_{n,p-n}\fred v_p$ according
to Lemma \ref{lem-red-prop}) and $x_{n+1,p-n}\fred^k v_{p+1}$ would
imply: $v_{p+1}\fred^{\ell-k} u$, contradicting the irreducibility
of $v_{p+1}$. Thus $\iota(x)\in\Sigma(\D)$.

\begin{figure}
$$
 \xymatrix{
     &   & v_p \\
     &  \ast \ar[ur]^{f^{\ell-k}} & w \ar[u]^i\\
   x_{n,p-n} \ar[ur]^{f^k} & v_{p+1} \ar[ur]^{f^{\ell-k}} \ar[u]^i \\
   x_{n+1,p-n} \ar[u]^i \ar[ur]^{f^k}
 }
 $$
\caption{Proof of $v_p\leadsto v_{p+1}$ for
Proposition \ref{prop-partial-iso}.}\label{fig-vpred}
\end{figure}

\medbreak

Let us prove that $\iota\circ\pi=\id_{\Sigma(\D)}$. Let
$v\in\Sigma(\D)$ and ${x}=\pi({v})\in\oX_V$.
Let us check that $x$ belongs to $\oX_V^M$. We
have, for $p\in\Z$ and $n\geq1$,
 $$
   x_{n,p-n}=(\pi(v_{p-n}v_{p-n+1}\dots))_n
    = f^{|v_{p-n}|}(y_{|v_{p-n}|+n})
 $$
where $y_{|v_{p-n}|+n}\fred^k v_p$ for $k=|v_{p-n}|+n-|v_p|$ by
Corollary \ref{coro-proj}. For $n\geq|v_p|$, $k\geq |v_{p-n}|$ and
 \begin{equation}\label{eq-x-red-v}
 x_{n,p-n}=f^{|v_{p-n}|}(y_{|v_{p-n}|+n})\fred v_p.
 \end{equation}
Thus ${x}$ is eventually Markov at any time $p$. ${x}\in\oX_V^M$
as claimed. Observe that eq. (\ref{eq-x-red-v}) also implies that
$\iota({x})={v}$, i.e., $\iota\circ\pi=\id_{\Sigma(\D)}$ as
claimed.

\medbreak

It remains to show that $\iota:\oX_V^M\to\Sigma(\D)$ is one-to-one.
Let ${x},{y}\in\oX_V^M$ with
$\iota({x})=\iota({y})=:{v}$. Let $p\in\Z$. As $v_p\arrow
v_{p+1}$, there is a unique $u_p^1$ such that $i(u^1_p)=v_p$ and
$u^1_p\fred^k v_{p+1}$ for $k=|v_p|+1-|v_{p+1}|$ by Lemma
\ref{lem-u-unique}.

For $n$ large enough, we have $x_{n,p-n},y_{n,p-n}\fred^\ell v_p$
for $\ell=n-|v_p|$ and $x_{n+1,p-n},y_{n+1,p-n}\fred v_{p+1}$.
Then $x_{n+1,p-n}\fred^\ell x_{n+1-\ell,p-n+\ell}=:w$ and $w$
must satisfy $i(w)=v_p$ and $w\fred v_{p+1}$ (observe that
$|w|=n+1-\ell= |v_p|+1\geq|v_{p+1}|$). By Lemma \ref{lem-u-unique},
$w=x_{n+1,p-n}$. Thus
 $$
    x_{n+1,p-n},y_{n+1,p-n}\fred
       u^1_p=x_{|v_p|+1,p-|v_p|}=y_{|v_p|+1,p-|v_p|}.
 $$
We want to repeat this analysis with $u^1$ replacing $v_p$.
First we check that $u^1_p\arrow u^1_{p+1}$, i.e., that there is some
$w$ such that $i(w)=u^1_p$ and $w\fred u^1_{p+1}$ (but some $u^1_p$
might be reducible). Indeed, $f^k:\itree(u^1_p)\to\itree(v_{p+1})$ is an
isomorphism so that there exists $w\in\itree(u^1_p)$ with
$f^k(w)=u^1_{p+1}$. Lemma \ref{lem-red-prop} gives then that
$w\fred u^1_{p+1}$. But this says that $u^1_p\arrow u^1_{p+1}$, as
claimed.

We proceed inductively, assuming that some sequence $(u^j_p)_{p\in\Z}$
has been obtained such that $u^j_p\leadsto u^j_{p+1}$ and
 \begin{equation}\label{eq:ujp}
    x_{n+j,p-n} \fred u^j_p = x_{|u^j_p|+1,p-|v_p|} = y_{|u^j_p|+1,p-|v_p|}
 \end{equation}
We define $u^{j+1}_p$ as the unique piece such that $i(u^{j+1}_p)=u^j_p$
and $u^{j+1}_p\fred u^j_{p+1}$. The same reasoning as above yields
(\ref{eq:ujp}) with $j$ replaced by $j+1$.

As $|u^j_p|=|v_p|+j\to\infty$, we obtain $x=y$.
\end{demo}

\begin{coro}\label{coro-h-pi}
The induced maps on the invariant probability measures
$\pi:\Prob(\sigma,\Sigma(\D))\to\Prob(F_V,\oX_V)$ and
$\pi:\Prob(\sigma,\Sigma_+(\D))\to\Prob(F_V,X_V)$ are one-to-one
and preserve ergodicity and entropy.
\end{coro}

\begin{demo}
That $\pi:\Sigma(\D)\to\oX_V$ is a partial isomorphism trivially
implies the stated properties of
$\pi:\Prob(\sigma,\Sigma(\D))\to\Prob(F_V,\oX_V)$. To finish,
recall that the natural extension construction preserves
ergodicity and entropy.
\end{demo}

%\begin{coro}\label{coro-countable-pi}
%The map $\pi:\Sigma_+(\mathcal D)\to X_V$ is at most
%countable-to-one almost everywhere w.r.t. any invariant
%probability measure on $\Sigma_+(\mathcal D)$.
%\end{coro}
%
%\begin{demo}
%By Poincar\'e recurrence,  for a.e. $v\in\Sigma_+(\mathcal D)$,
%$L_v:=\liminf_{p\to\infty}|v_p|<\infty$. Hence, there exist
%$p_v\geq L_v$ such that $|v_{p_v}|=L_v$. Now, for $p\geq p_v$,
%$|v_p|\leq L+(p-p_v)$ so (\ref{eq-iso}) implies that, writing $x:=
%\pi(v)$, $v_p=f^k(x_{|v_p|,p-|v_p|})$ and $p-|v_p|\geq
%p_v-L\geq0$. Therefore $x=\pi(v)$, $p_v$ and $v_0,\dots,v_{p_v}$
%determine $v$. As there are countably many choices of the finitely
%many vertices $v_0,\dots,v_{p_v}$, the corollary is proved.
%\end{demo}

\subsection{Entropy of the non-Markov part}

\begin{prop} \label{prop-non-M}
If $\mu$ is an invariant and probability measure with
$\mu(\oX_V\setminus \oX_V^M)=0$, then
 $$
    h(F_V,\mu)\leq h_\C(V).
 $$
\end{prop}

To analyze the non-Markov part, the first step is the following:

\begin{lem}\label{lem-non-M}
Up to a set of zero measure with respect to any invariant
probability measure, each ${x}\in\oX_V\setminus\oX_V^M$ satisfies:
for {\bf all} $p\in\Z$ there exist arbitrarily large integers $n$
such that $x_{n,p-n}$ is an $f$-irreducible vertex.
\end{lem}

\begin{demof}{the Lemma}
By definition ${x}\in\oX_V\setminus\oX_V^M$ iff there {\em exists}
$p\in\Z$ as in the statement of the Lemma. Let $X(p)$ be the set
of such ${x}$. The lemma is clearly equivalent to the fact that,
for any invariant probability $\mu$,
 \begin{equation}\label{eq-union-inter}
   \mu\left(\bigcup_{p\in\Z}X(p)\right) =
   \mu\left(\bigcap_{p\in\Z}X(p)\right).
 \end{equation}
It is enough to prove this for ergodic $\mu$'s such that the union
has positive and hence full measure. If we prove that
$X(p+1)\subset X(p)$, it will follow that
$\mu\left(\bigcup_{p\in\Z}X(p)\right)=\lim_{p\to-\infty} \mu(X(p))$
which is equal to $\mu(X(p))$ for any $p\in\Z$ by invariance of
$\mu$, proving (\ref{eq-union-inter}). But observe that by Lemma
\ref{lem-red-prop}
 $$
   x_{n,p-n} \fred x_{N,p-N}
    \implies x_{n+1,p+1-(n+1)} \fred x_{N+1,p+1-(N+1)})
 $$
so that $x\notin X(p)\implies x\notin X(p+1)$, which concludes the
proof.
\end{demof}

Recall that the entropy of an invariant and ergodic probability measure
$\mu$ can be computed as \cite{Katok}:
 \begin{multline*}
   h(F_V,\mu) = \lim_{\eps\to0} h(F_V,\mu,\eps)
   \text{ with }\\
   h(F_V,\mu,\eps) = \limsup_{n\to\infty} \frac1n\log\min\left\{\# S:
     \mu\bigl( \bigcup_{x\in S} B(x,\eps,n) \bigr)>\mu_0 \right\}
 \end{multline*}
where $0<\mu_0<1$ is arbitrary.

\begin{demof}{the Proposition}
Let $\mu$ be an invariant probability measure carried by
$\oX_V\setminus\oX_V^M$. We may and do assume that $\mu$ is
ergodic. Let $\alpha>0$ be some small number. There exists $r>0$
(depending on $\mu$) such that $h(F_V,\mu)\leq h(F_V,\mu,r)+\alpha$.
Fix $L_1<\infty$ and $r_1>0$ such that for ${x},{y}\in\oX_V$,
$d(x_{2L_1,-L_1},y_{2L_1,-L_1})<r_1 \implies d({x},{y})<r$ (for
any distance on $X_V$ and $\oX_V$ compatible with the
topologies). Let $L_2$ be such that $r(r_1,n,\C_n)\leq
e^{(h_\C(V)+\alpha)n}$ for all $n\geq L_2$ and fix, for each such
$n$, some $(r_1,n)$-cover\footnote{Recall the definition
of the entropy of sequences.} $C_n$ of $\C_n$
with this minimum cardinality. For each $v\in C_n$, we pick some
$x\in \oX_V$ such that $x_{n,0}=v$ and let $X_n:=\{x^v\in
X_V: v\in C_n\}$.

Let $L>>L_1\log K/2\alpha+L_2$ where $K$ is the minimum cardinality
of an $r$-dense subset of $\oX_V$. It follows from Lemma
\ref{lem-non-M} that there exists a measurable function
$n:\oX_V\to\N$ such that, for $\mu$-a.e. $\un{x}\in\oX_V$:
 \begin{itemize}
  \item $n(\un{x})\geq L$;
  \item $x_{n,-n}$ is $f$-irreducible for $n=n(\un{x})$.
 \end{itemize}
Hence (see \cite[p. 394]{QFT}) a $\mu$-typical $x$ satisfies the
following. For all large $n$, there exist disjoint integer
intervals $[a_i,b_i)\subset[0,n)$, $i=1,\dots,s$, such that
 \begin{enumerate}
  \item $\sum_{i=1}^s b_i-a_i\geq (1-\alpha)n$;
  \item $b_i-a_i\geq L$ for all $i=1,\dots,s$;
  \item $x_{b_i-a_i+1,a_i}\not\fred x_{b_i-a_i,a_i+1}$: in particular,
  $x_{b_i-a_i+1,a_i}$ is $f$-irreducible. Thus
  $F_V^{a_i+L_1}({x})\in B(y,r,b_i-a_i-2L_1)$ for some $y\in X_{b_i-a_i+1}$.
 \end{enumerate}
It follows (see, e.g., the same reference)
that $h(F_V,\mu)\leq h_\C(V)+3\alpha+\alpha|\log\alpha|\leq H$. As
$\alpha>0$ is arbitrarily small, this concludes the proof.
\end{demof}

\subsection{Entropy at infinity in the diagram}

\begin{prop} \label{prop-h-at-infty}
Let $H>h_C(V)+h_\wloc(V)$. Then there exists a finite subset
$\D_0\subset\D$ such that:
 $$
   h(\D\setminus\D_0,\D):=\inf_{\mu_0>0} \sup \left\{ h(\sigma,\mu):
     \mu\in\Prob_\erg(\Sigma(\D)) \text{ and } \mu\left(\bigcup_{D\in\D_0}[D]\right)<\mu_0
     \right\} \leq H
 $$
where $\Prob_\erg(\Sigma(\D))$ is the set of shift-invariant and
ergodic probability measures on $\Sigma(\D)$.
\end{prop}

\begin{demo}
It is enough to find $\D_0$ and $\mu_0>0$ such that if
$\mu\in\Prob_\erg(\Sigma(\D))$ satisfies:
 \begin{equation}\label{eq-at-infty}
 \mu\left(\bigcup_{D\in\D_0}[D]\right)<\mu_0,
 \end{equation}
then $h(\sigma,\mu)\leq H$.

Let $\alpha>0$ be so small that $h_\C(V)+ h_\wloc(V)+ 4\alpha+
\alpha|\log\alpha|\leq H$. Let $r>0$ be such that, for all invariant and
ergodic probability measures $\mu$ with $h(F_V,\mu)>h_\C(V)$:
 $$h(F_v,\mu)-h(F_V,\mu,r)\leq h_\wloc(F_V)+\alpha$$
(the point here is that $r$ and therefore $L_1$ and $r_1$ are now fixed,
especially they are independent from $\mu$
---compare with Proposition \ref{prop-partial-iso}). Fix
$L_1<\infty$ and $r_1>0$ such that, for all $x,y\in\oX_V$,
$d(x_{2L_1,-L_1},y_{2L_1,-L_1})<r_1 \implies d({x},{y})<r$.
We increase $L_1$ if necessary so that $L_1>r_1^{-1}$. Recall that
$\C=(C_n)_{n\geq1}$ with $C_n$ the set of irreducible vertices of
order $n$. Let $L_2$ such that $r(r_1,n,\C_n)\leq
e^{(h_\top(\C)+\alpha)n}$ for all $n\geq L_2$. Let $K$ be the
cardinality of a finite $r$-dense subset of $\oX_V$ and let
$L>\alpha^{-1}L_1\log K+L_2$.

Finally let
 $$
    \D_0=\{v\in\D:|v|\leq L\}
 $$
and let $\mu_0>0$ be a very small number to be specified later.

Let $\mu$ be an ergodic invariant probability measure on
$\Sigma(\D)$ satisfying (\ref{eq-at-infty}). We bound
$h(\sigma,\mu)$. First observe that by Corollary \ref{coro-h-pi},
$h(\sigma,\mu)=h(F_V,\pi_*\mu)$. Let ${x}\in\oX_V$ be a
$\pi_*\mu$-typical point. Thus ${x}=\pi({v})$ with ${v}$ a path on
$\D$ spending a fraction of its time less than $\mu_0$ in $\D_0$.

This implies that there exist disjoint integer intervals
$[a_1,b_1),\dots\subset[0,n)$ with $v_{b_i}
\in\D\setminus\D_0$ such that $\sum_i b_i-a_i\geq(1-\mu_0)n$
and $x_{|v_{b_i}|,b_i-|v_{b_i}|}=v_{b_i}$. The latter implies:
$d(x_{a_i+k},f^k(v_{b_i}))<r_1$ for all $k\in[0,b_i-a_i-L_1)$. Note
that the $b_i-a_i$ are large (larger than $L$). By definition of
$\D$, the $v_{b_i}$s are $f$-irreducible. It follows as in the proof
of Proposition \ref{prop-non-M} that
 $$
   h(F_V,\pi_*\mu,r)\leq h_\C(V)+\alpha+\frac1L+\frac1L|\log \frac1L|+(\mu_0+2L_1/L)\log K
     \leq h_\C(V)+3\alpha+\alpha|\log\alpha|
 $$
if $\mu_0=\mu_0(V,r,\alpha)$ is small enough. If $h(F_V,\pi_*\mu)\leq h_\C(V)\leq H$,
there is nothing to show. Otherwise,
 $$
   h(F_V,\pi_*\mu)\leq h_\C(V)+h_\wloc(F_V,r) +4\alpha+\alpha|\log\alpha|\leq H
 $$
as claimed.
\end{demo}

\subsection{Conclusion of the Analysis of Large Entropy Measures}

We collect all the partial results and check that they imply the
first two claims of Theorem \ref{theo-struct}.

First, let $V$ be a \sQFT puzzle. Propositions
\ref{prop-partial-iso} and \ref{prop-non-M} immediately imply that
$\oX_V$ is entropy-conjugate with constant $h_\C(V)$ to the Markov
shift, $\Sigma(\D)$, proving claim (1) of the Theorem.

For claim (2), we assume that $V$ is \QFT:
$h_\top(X_V)>H_*:=h_\C(V)+h_\wloc(X_V)$.  Proposition
\ref{prop-h-at-infty} implies that $h_\infty(\D)\leq H_*$.

Take $H$ strictly between $H_*$ and $h_\top(V)$: $\D$ contains
only finitely many irreducible Markov subshifts $S$ with entropy
$h(S)\geq H$. This implies that $h(S)>H_*\geq h_\infty(\D)\geq
h_\infty(S)$. Hence, by the result of Gurevi\v{c} and Zargaryan
\cite{GZ} quoted in Proposition \ref{prop-ruette} below these
irreducible subshifts are SPR. This proves claim (2) of Theorem
\ref{theo-struct}.

%%%%%%%%%%%%%%%%%%%%%%%%%%%%%%%%%%%%%%%%%%%%%%%%%%%%%%%%%%%%%%%%%%%%%%%%%%%%%%%%%%
\section{Periodic Structure}\label{sec-period}
%%%%%%%%%%%%%%%%%%%%%%%%%%%%%%%%%%%%%%%%%%%%%%%%%%%%%%%%%%%%%%%%%%%%%%%%%%%%%%%%%%

In this section we prove Claim (3) of Theorem \ref{theo-struct}
which relates most periodic orbits in the Markov shift with most
periodic orbits in some fine scale approximation $i_N(X_V)$ of the
puzzle dynamics $X_V$. It is here that we need determinacy,
exactly once, to prove eq. (\ref{eq-fix1}).

\subsection{Partition of the periodic points}

The proof will use two integer parameters $N,L\geq1$
depending on $\eps>0$. We shall denote $i_N\circ\pi:\Sigma(\D)\to
i_N(\oX_V)$ by $\pi_N$. The $n$-fixed points $\xi=\sigma^n(\xi)$ of
$i_N(X_V)$ satisfy exactly one of the following properties:
\begin{enumerate}
 \item[(P1)] there exist $v\in\pi_N^{-1}(\xi)\in\Sigma(\D)$
 such that $I_N(v):=\{p\geq0:|v_p|<N\}$ is infinite.
 \item[(P2)] $\pi_N^{-1}(\xi)\ne\emptyset$ but for
 all $v$ in this set,  $I_N(v)$ is finite.
 \item[(P3)] $\pi_N^{-1}(\xi)=\emptyset$.
\end{enumerate}
Denote by $\tFix_i(n)$, $i=1,2,3$, the corresponding sets of
periodic points of $i_N(X_V)$ (these sets do not
depend on $L$, which will define a splitting of $\tFix_2(n)$ below).

On the other hand, we consider on the Markov shift only the
periodic points defined by \new{low loops}:
 $$
   \widehat\Fix_1(n) := \{ v\in\Sigma(\D): \sigma^nv=v
      \text{ and } \{v_0,\dots,v_{n-1}\}\cap \D_N\ne\emptyset\}.
 $$
We shall say nothing about the others.

\subsection{Low loops and periodic points of $i_N(X_V)$}

Let $\eps>0$, $N_0$ and $\D_0$ be given as in the statement of the
Theorem. Fix $N\geq N_0$ so that $\D_N\supset\D_0$ and
$h(\D\setminus\D_N)\leq h_\C(V)+h_\wloc(V)+\eps/2$ (which is
possible by Proposition \ref{prop-h-at-infty} as
$h(\D\setminus\D_N)\leq h(\D\setminus\D_N,\D)$).

We first claim that for all $n\geq1$:
 \begin{equation}\label{eq-fix1}
     \#\tFix_1(n)=\#\widehat\Fix_1(n).
 \end{equation}

We need the following consequence of determinacy:

\begin{lem}\label{lem-determined}
Let $V$ be a determined puzzle and $N\geq1$. Let
$v,v'\in\Sigma(\mathcal D)$. If $x=\pi(v)$ and $x'=\pi(v')$
satisfy $i_1(x)=i_1(x')$, then:
 \begin{equation}\label{eq-x-det}
      v_0=v'_0\implies \forall n\geq0\; x_{n,-n} = x'_{n,-n}
 \end{equation}
\end{lem}

\begin{demo}
For $n=n_0:=|v_0|$, the right hand side of (\ref{eq-x-det})
follows from $v_0=x_{n_0,-n_0}= x'_{n_0,-n_0}$, which holds by
(\ref{eq-iso}). This implies (\ref{eq-x-det}) for $n\leq n_0$.
Assuming it for some $n\geq n_0$, (\ref{eq-iso}) again implies
$x_{n+1,-n-1},x'_{n+1,-n-1} \fred x_{n,-n} =x'_{n,-n}$. Together
with the determinacy and $i_1(x_{n+1,-n-1})=i_1(x'_{n+1,-n-1})$,
this completes the induction and the proof of the lemma.
\end{demo}

We deduce (\ref{eq-fix1}) from this Lemma. Let
 $$
   \Sigma(N):=\{ v\in \Sigma(\D): \exists p\to\infty\; |v_{p}|<N \}.
 $$
By Lemma \ref{lem-determined}, $\pi_N|\Sigma(N)$ is one-to-one.
$\tFix_1(n)$ is by definition the set of fixed points $\xi$ of
$\sigma^n$ in $\pi_N(\Sigma(N))$. By the injectivity of
$\pi_N|\Sigma(N)$ and the $\sigma$-invariance of $\Sigma(N)$, such
$\xi$ are the $\pi_N$ images of the fixed points of $\sigma^n$ in
$\pi_N(\Sigma(N))$. This proves the claim  (\ref{eq-fix1}).

\medbreak

\subsection{Remaining loops and periodic points}

To conclude we check that the remaining factors are holomorphic.

\begin{lem}\label{lem-zeta23}
For every $\eps>0$, there exists $N_0(\eps)<\infty$ such that, for all
$N\geq N_0(\eps)$, for all $n\geq 1$:
 \begin{gather}
  \label{eq-zeta2}
    \#\tFix_2(n) \leq C(N,\eps)
  e^{(h_\C(V)+\eps)n}\\
  \label{eq-zeta3}    \#\tFix_3(n) \leq C(N,\eps)
  e^{(h_\C(V)+\eps)n}.
\end{gather}
\end{lem}

\begin{remk}
Of course, $N$ being large, $h(\Sigma(\D\setminus\D_N))<H:=
h_\C(V)+ h_\wloc(V)+ \eps <h_\top(V)$ but this is unsufficient to
prove (\ref{eq-zeta2}) as Gurevi\v{c} entropy only controls the
number of loops based at a \emph{fixed} vertex. Indeed, in some
examples, $\D\setminus\D_N$ contains infinitely many loops of each
length.
\end{remk}

This Lemma will be enough to conclude the proof of Theorem
\ref{theo-qft-zeta}.

\subsection{Paths above $N$}
To prove (\ref{eq-zeta2}) we bound the number of the $n$-periodic
projections to $i_N(X_V)$ of (not necessarily periodic) paths on
$\D\setminus\D_N$.

The proof is similar to that of Proposition \ref{prop-h-at-infty}.
For $n\geq N$, let
  \begin{multline}\label{eq-C-n}
    C(n):=\biggl\{ (i_N(v),i_N(f(v)),\dots,i_N(f^{n-N}(v))) \in V_{\leq
    N}^{n-N}:
       v\in \C_n  \biggr\}.
  \end{multline}
By (\ref{eq-iso}), for all $w\in\Sigma(\D)$, writing
$y:=\pi(w)\in\oX_V$,
 \begin{equation}\label{eq-shadow}
      w_0 \text{ irreducible with }\ell:=|w_0|\geq N
         \implies (y_{N,-\ell+1},y_{N,-\ell+2},\dots,y_{N,-N})\in
         C(\ell)
 \end{equation}
so that $\{C(n):n\geq 1\}$ controls the projections of high paths.
By definition,
 \begin{equation}\label{eq:C-bound}
     \# C(n)\leq C(N,\eps) e^{(h_\C(V)+\eps)n},
 \end{equation}
hence this control should give an
entropy bound. Let us see the details (note that the bound proved below is
only $h_\C(V)+h_\wloc(V)$). We fix $L=L(N,\eps)$, a large integer.

Let $n\geq1$ and $\xi\in i_N(\oX_V)$ with $\sigma^n(\xi)=\xi$
satisfying (P2): $\xi=\pi_N(v)$ for some $v\in\Sigma(\D)$
satisfying: $|v_p|\geq N$ for all $p\geq p_0$. By periodicity of
$\xi$, we can assume $p_0=0$ by shifting $v$ a multiple of the
period. We shift again $\xi$ to ensure $I_L(v)=\emptyset$ if
$I_L(v)$ is finite, $|v_0|\leq L$ otherwise (this might produce an
irrelevant factor $n$ in the estimates). Let $x=\pi_N(v)\in\oX_V$.

Define inductively the integers $r\geq1$, $n>b_1>\dots>b_r\geq 0$
as follows. $b_1:=n-1$. If $b_i-|v_{b_i}|\geq 0$ and
$|v_{b_i}|>L$, then set $b_{i+1}:=b_i-|v_{b_i}|\geq0$, otherwise let
$r:=i$. Finally set $\ell_i:=|v_{b_i}|$ for all $i=1,2,\dots,r$.

Now, by (\ref{eq-shadow}), $x_{N,b_i-\ell_i+k}=i_N \circ
f^k(v_{b_i})$ for $0\leq k<\ell_i-N$.  Thus, there exists $x(i)\in
C(\ell_i)$, such that $\xi_{b_i-\ell_i+k}=(x(i))_k$ for $0\leq
k\leq\ell_i-N$.

Notice that $\ell_i\geq L$ for $i=1,\dots,r-1$. Hence, given
$n\geq1$ and $b_r$, when $v$ ranges over $\Sigma(\D\setminus
\D_N)$, the number of choices for the integers $b_1,\dots,b_{r-1}$
is at most $e^{\eps(n-b_r)}$ as $L$ is large. On the other hand,
eq. (\ref{eq:C-bound}) yields, for any $\ell\geq0$,
 $$
    \#\{\xi_{-\ell+1}\dots\xi_{0}:\xi=\pi_N(v) \text{ with } v\in\Sigma(\D)
    \text{ and }
     |v_0|=\ell\} \leq C(N,\eps) e^{(h_{\C}(V)+\eps)\ell} \times (\# V_{\leq N})^N
  $$
 This implies:
 \begin{multline}\label{eq-case1}
    \#\{\xi_{b_r}\dots\xi_{n-1}:\xi\in\tFix_2(n) \text{ with given $b_r$ and }
     I_N(v)=\emptyset\}
     \leq \\
     (\# V_{\leq N})^{N(n-b_r)/L} \cdot e^{\eps (n-b_r)}\cdot C(N,\eps)^{(n-b_r)/L} \exp
       (h_\C(V)+\eps)(n-b_r) \\
          \leq C e^{(h_C(V)+3\eps)(n-b_r)}
 \end{multline}
using that $L=L(N,\eps)$ is large. Note that there are at most $n$
possibilities for $b_r$. It remains to count the possibilities for
$\xi_0\dots\xi_{b_r-1}$. There are two cases.

First case: $|v_{b_r}|\leq L$.
$\tFix_2^{(1)}(n)$ be the corresponding subset of $\tFix_2(n)$.
Recall that in this case $I_L(v)$ is infinite and $|v_0|\leq L$.
Hence $v_0\dots v_{b_r}$ is a path on $\D\setminus\D_N$ which starts and
ends in the finite subgraph $\D_L\setminus\D_N$.  The number of such paths with given $b_r$  is bounded using the Gurevi\v{c} entropy by:
 $$
    (\#\D_L\setminus\D_N)^2     e^{(h(\D\setminus\D_N)+\eps/2)b_r}
      \leq (\#\D_L\setminus\D_N)^2 e^{(h(V)+h_\wloc(V)+\eps)b_r}.
  $$
Combining with eq. (\ref{eq-case1}) and summing over $1\leq r\leq n/L+1$ and $0\leq b_r\leq n$ we obtain:
 $$
    \#\tFix_2^{(1)}(n) \leq C'(L,N,\eps) n^2 e^{(h(V)+h_\wloc(V)+3\eps)n}.
  $$

Second case: $|v_{b_r}|>L$ and $b_r-|v_{b_r}|=:-\delta<0$. Let
$\tFix_2^{(2)}(n)$ be the corresponding subset of $\tFix_2(n)$. We
shift $\xi$ by $\delta$ (and add $\delta$ to each $b_i$) so
$b_r-|v_{b_r}|=0$ (doing this we lose the property $|v_0|\leq L$
if $I_L(v)$ is infinite). Of
course, $b_1,\dots,b_s\geq n$ for some (maximum) $s\geq1$. We
forget about $b_1,\dots,b_{s-1}$ and we trim $v_{b_s}$ in the
following way. We replace $b_s$, $v_{b_s}$ and $\ell_s$ by $n-1$,
$v_*:=i^{b_s-n+1}(v_{b_s})$ and $\ell_*:=\ell_s-(b_s-n+1)\geq1$. We
have now that $x_{-N,b_s-\ell_s+i}=i_N(f^i v_*)$ for $0\leq
i<\ell_*-N$. We may now apply (\ref{eq-case1}) with $b_r=0$.
This concludes the proof of (\ref{eq-zeta2}).

\subsection{Unliftable periodic orbits}

Let $\xi\in\tFix_3(n)$. By definition, for any $x\in\oX_V$
projecting to $\xi$, there exist $p\in\Z$ and arbitrarily large
integers $k$ such that, $x_{k,p-k}$ is $f$-irreducible. Take such
an integer $k\geq N+n$ and observe that $v_*:=x_{N+n,p-k}$ is
$f$-irreducible by Lemma \ref{lem-red-prop}. Hence, setting
$q:=p-k-|v_*|$,
 $$
   \xi_{q+i} = i_N(f^i(v_*)) \text{ for all }0\leq i<n.
 $$
This implies that
 $$
   \#\tFix_3(n) \leq C(N,\eps) e^{(h_\C(V)+\eps)n},
 $$
proving eq. (\ref{eq-zeta3}) and concluding the proof of Theorem
\ref{theo-qft-zeta}.

\section{Semi-local zeta functions of SPR Markov
shifts}\label{sec-Markov-zeta}\label{sec:semilocal-zeta}

We give a proof of Theorem \ref{theo-zeta} about the meromorphy of
the semi-local zeta functions of SPR Markov shifts, after
recalling the relation between the entropy at infinity and the SPR
property.

\subsection{SPR property and entropy at infinity}

A combinatorial quantity related to our $h_\infty(G)$ first
appeared in a work of Gurevi\v{c}-Zargaryan \cite{GZ} to give
a sufficient condition for being SPR, which was then shown to be
necessary (see \cite[Theorem 3.8]{GS}). The explicit relation
(\ref{eq-ruette}) below is due to Ruette \cite{RuetteClass}.

\begin{prop}[Gurevi\v{c}-Zargaryan, Gurevi\v{c}-Savchenko, Ruette]\label{prop-ruette}
Let $G$ be a countable, oriented, irreducible graph with
$h(G)<\infty$. The graph $G$ is SPR iff $h_\infty(G)<h(G)$ where
the entropy at infinity $h_\infty(G)$ has been defined in
\ref{eq-h-infty} and can be computed as:
 \begin{multline}\label{eq-ruette}
    h_\infty(G) = \inf_{F\subset\subset G} \max_{u,v\in F}
    \limsup_{n\to\infty}\frac1n\log\#\{(x_0,\dots,x_{n})\in \{u\}\times G^{n-1}
      \times \{v\}: \\
      \forall i=0,\dots,n\; x_i\to x_{i+1}\text{ on }G
      \text{ and }
      \{x_1,x_2,\dots,x_{n-1}\}\cap F =\emptyset\}
 \end{multline}
where $F\subset\subset G$ means that $F$ ranges over the finite
subgraphs of $G$.
\end{prop}

Observe that by this proposition,
the conclusion of our Theorem \ref{theo-zeta} is
non-trivial iff the Markov shift is SPR.

\subsection{Semi-local zeta function of large subsets}
The first step of the proof of Theorem \ref{theo-zeta} is the
following weaker claim on semi-local zeta functions defined by
large subgraphs:

\begin{claim}\label{claim-large-zeta}
For every $\eps>0$, there exists a finite subset
$F_0\subset\subset G$ such that for all finite subsets $F_0\subset
F\subset\subset G$, the semi-local zeta function $\zeta^G_F(z)$
extends meromorphically to $|z|<\exp-(h_\infty(G)+\eps)$.
\end{claim}

The crux of the proof is to check that
 $$
   \zeta^G_F(z) = 1/\det(\operatorname{Id}-L(z))
 $$
with $L(z)$ a finite square matrix with holomorphic entries for
$|z|<e^{-(h_\infty(G)+\eps)}$ and $\operatorname{Id}$ is the identity matrix.

One can give a direct, self-contained proof by
generalizing an algebraic formula for the determinant of a finite
matrix in terms of its block decomposition used for a similar purpose in
\cite{hofbauer-keller}. We give a shorter proof based on the formula in
eq. (\ref{eq:BGMY}), quoted from \cite{BGMY}, as suggested by a referee.

\begin{demof}{Claim \ref{claim-large-zeta}}
Let $F$ be a finite subgraph of $G$. For each $u,v\in F$, define
 \begin{multline*}
    f^F_n(u,v):= \#\{x_1\dots
   x_{n-1}:x\in\Sigma(G) \text{ s.t. }
      x_0=u,\; x_n=v\\ \text{ and }
      \{x_1,x_2,\dots,x_{n-1}\}\cap F=\emptyset\}.
 \end{multline*}
Recall the definition of $h(G\setminus F,G)$ as introduced in Proposition \ref{prop-h-at-infty}.
There is $F_0\subset\subset G$ such that, if $F_0\subset F\subset\subset G$,
then $h(G\setminus F,G)<h_\infty(G)+\eps$ for arbitrarily small $\eps>0$. Hence,
by eq. (\ref{eq-ruette}), for all
$u,v\in F$:
 \begin{equation}\label{eq:fF}
   \limsup_{n\to\infty}\frac1n\log f^F_n(u,v)<
      h_\infty(G)+\eps.
 \end{equation}

Now define $L_n=L_n(z)$ to be the  $F\times F$-matrix
with following polynomial entries in $z$:
 $$
    L_n(u,v) := \sum_{0\leq k\leq n} f^F_k(u,v)z^k.
 $$
Consider the zeta function:
 \begin{multline*}
   \zeta^F_n(z):=\exp \sum_{p\geq1} \frac{z^p}p \#\{x_0\dots x_p:
     x\in\Sigma(G),\; \sigma^p(x)=x,\\
    \{i\in\Z:x_i\in F\}
     \text{ has gaps of lengths at most }n\}.
 \end{multline*}
The formula from \cite{BGMY} is:
 \begin{equation}\label{eq:BGMY}
   \zeta^F_n(z)=\det(\operatorname{Id}-L_n(z))^{-1}.
 \end{equation}

Let $L$ be the $F\times F$ matrix with power series entries defined by
$L(u,v):=\lim_{n\to\infty} L_n(u,v)$. The limit here is in the sense of formal power
series. Increasing $n$ only adds high powers, hence this limit exists. Also $\zeta^G_F(z)=\lim_{n\to\infty} \zeta^F_n(z)$.
By eq. (\ref{eq:BGMY}), we get $\zeta^G_F(z)=\det(\operatorname{Id}-L(z))^{-1}$.

Eq. (\ref{eq:fF})implies that the entries of $L$ have a radius of
convergence at least $e^{-h_\infty(G)-\eps}$. Therefore $\zeta^G_F(z)$
must be meromorphic over $|z|<\exp -(h_\infty(G)+\eps)$.
\end{demof}

\subsection{Proof of Theorem  \ref{theo-zeta}}
Claim \ref{claim-large-zeta} shows the meromorphy of semi-local
zeta functions relative to large finite subsets.

We first show the last claim of the Theorem. For $\eps>0$, let
$F_0\subset\subset G$ as in the proof of Claim
\ref{claim-large-zeta}. let $F',F\supset F_0$ be other finite
subgraphs. Let $F_1:=F\cup F'$.
 $$
    \frac{\zeta^G_{F_1 }(z)}{\zeta^G_{F}(z)}
     = \exp \sum_{n\geq1} \frac{z^n}{n}
    \#\{x\in\Sigma:\sigma^n(x)=x \text{ and }
    \{x_0,\dots,x_{n-1}\}\text{meets $F_1$ but not $F$} \}.
 $$
The radius of convergence of the above series is at least
$e^{-h(G\setminus F_0)}\geq e^{-(h_\infty(G)+\eps)}$. The same
applies to $\zeta^G_{F_1 }(z)/\zeta^G_{F'}(z)$. This proves that
$\zeta^G_F(z)/\zeta^G_{F'}(z)$ is a holomorphic non-zero function
over $|z|<e^{-(h_\infty(G)+\eps)}$.

We now show that all semi-local zeta functions are meromorphic on
$|z|<e^{-h_\infty(G)}$, finishing the proof of Theorem
\ref{theo-zeta}.

Let $F\subset\subset G$. Let $\eps>0$. By taking $H$ such that
$F\subset H\subset\subset G$, $H$ large enough, we can ensure that
 $$
    h(G\setminus H,G\setminus F) \leq h_\infty(G\setminus F) + \eps
 $$
(see the definition of $h(\cdot,\cdot)$ in Proposition  \ref{prop-h-at-infty}). After
possibly increasing $H$, Claim \ref{claim-large-zeta} ensures that
$\zeta^{G\setminus F}_{H\setminus F}(z)$ has a meromorphic
extension to $|z|<\exp-h_\infty(G\setminus F)-\eps$. We compute:
 \begin{equation}\label{eq-ratio-zeta}
  \begin{aligned}
     \zeta^G_H(z)/\zeta^G_F(z) &= \exp \sum_{n\geq 1} \frac{z^n}n \#\{x\in\Sigma(G):
        \sigma^nx=x \text{ and } \\
         &\qquad\qquad\qquad\qquad \{x_0,\dots,x_{n-1}\}\cap H\ne\emptyset \text{ but }
          \{x_0,\dots,x_{n-1}\}\cap F=\emptyset \} \\
        &= \exp \sum_{n\geq 1} \frac{z^n}n \#\{x\in\Sigma(G\setminus F):
        \sigma^nx=x \text{ and } \\
        &\qquad\qquad\qquad\qquad \{x_0,\dots,x_{n-1}\}\cap (H\setminus
        F)\ne\emptyset \} \\
        &= \zeta^{G\setminus F}_{H\setminus F}(z).
 \end{aligned}
 \end{equation}
As $h_\infty(G\setminus H)\leq h_\infty(G\setminus F)\leq
h_\infty(G)$ (see the remark after Definition \ref{defi-h-infty}).
Thus $\zeta_F^G(z)=\zeta_H^G(z)/\zeta^{G\setminus F}_{H\setminus
F}(z)$ is meromorphic on $|z|<\exp - h_\infty(G)-\eps$.
Letting $\eps>0$ decrease to $0$, finish the proof of Theorem
\ref{theo-zeta}.

\section{Proof of the Consequences}\label{sec-coro}

\subsection{Measures of maximum entropy}\label{sec-max-meas}
Our Structure Theorem implies that the set of maximum measures for
a \sQFT puzzle or for the associated Markov shift have the same
cardinality. We apply some results of Gurevi\v{c}.

First, according to \cite{Gurevic}, each irreducible subshift of a
Markov shift carries at most one maximum measure and this measure,
if it exists, is a Markov measure (which implies by
\cite{Bernoulli} that it is a finite extension of a Bernoulli).
Hence, a \sQFT puzzle has at most countably many maximum measures
(because it has at most countably many states) and a \QFT puzzle
has only finitely many irreducible components (because its
spectral decomposition contains finitely many irreducible
subshifts with maximum entropy).

The existence of a maximum measure for a \QFT puzzle follows from
the fact that the spectral decomposition of its Markov diagram must contain an
irreducible subshift with entropy equal to that of the puzzle and
that this subshift must be SPR. SPR implies positive recurrence which
is equivalent to the existence of a maximum measure by the same
result of Gurevi\v{c}. Theorem \ref{theo-maxmes} is proved.

\subsection{Zeta functions}\label{sec-qft-zeta}
We prove Theorem \ref{theo-qft-zeta}.
Recall that for the results involving the counting of the periodic
points, we assume, in addition to \QFT, determinacy. For
simplicity, we assume that the Markov diagram $\D$ is irreducible
and leave the general case to the reader. Let $\eps>0$. Theorem
\ref{theo-struct} gives a large integer $N$ such that the
$n$-periodic orbits of $i_N(X_V)$ and the loops of $\D$ going
through $\D_N$ can be identified up to an error bounded by $\exp
(h_\C(V)+ h_\wloc(V)+\eps)n$. Hence
 $$
    \zeta^N(z) := \exp \sum_{n\geq1} \frac{z^n}n \#\{x\in
    i_N(X_V):F_V^n(x)=x\}
 $$
is equal to the semi-local zeta function of $\D$ at $\D_N$ up to a
holomorphic, non-zero factor on the disk
$|z|<e^{-h_\C(V)-h_\wloc(V)-\eps}$ by Claim (3) of Theorem
\ref{theo-struct}. By Theorem \ref{theo-zeta}, this semi-local
zeta function can be extended to a meromorphic function on
$|z|<e^{-h_\C(V)-h_\wloc(V)}$, proving the main claim.

The singularities of $\zeta^N(z)$ on $|z|=e^{-h_\top(V)}$ are as claimed by
the same statement proved for local zeta function ($F$ reduced to
one vertex) by Gurevic and Savchenko \cite{GS}.

This concludes the proof of Theorem \ref{theo-qft-zeta}.

\subsection{Equidistribution of periodic points}

We give a sketch of the proof of Theorem \ref{theo-equi-per}
which is essentially that from
\cite{QFT} using the estimates of the analysis of the zeta
function above.

There is equidistribution for an irreducible SPR Markov shift
according to Gurevic and Savchenko \cite{GS}. For
the (easy) extension to the general case, it is enough to see
(like in \cite{QFT}) that the number of $n$-periodic points living
on an irreducible SPR Markov shift $\Sigma$ with period $p$ is
equivalent to $pe^{n h_\top(\Sigma)}$ if $n$ is a multiple of $p$,
zero otherwise.

To apply it to the puzzle, one has to recall the following facts
from the above analysis of the zeta function:
 \begin{itemize}
  \item the projection $\Sigma(\D)\to\oX_V$ is continuous;
  \item there is a one-to-one, period-preserving correspondence between
   $i_N$-projections of periodic points going through a large finite subset $F$ and a
   subset of the periodic points of $X_V$;
  \item the remaining periodic points both on $i_N(X_V)$ and $\Sigma(\D)$ contributes negligibly
   to the considered measures by the reasoning in the proof of
   Theorem \ref{theo-qft-zeta}.
  \end{itemize}

%%%%%%%%%%%%%%%%%%%%%%%%%%%%%%%%%%%%%%%%%%%%%%%%%%%%%%%%%%%%%%%%%%%%%%%%%%%%%
\section{Application to entropy-expanding maps}\label{sec-hexp}
%%%%%%%%%%%%%%%%%%%%%%%%%%%%%%%%%%%%%%%%%%%%%%%%%%%%%%%%%%%%%%%%%%%%%%%%%%%%%

We prove Theorem \ref{theo-hexp}: smooth entropy-expanding
maps introduced in \cite{BullSMF} define determined puzzles of
quasi-finite type, provided that they are endowed with a good
partition in the sense of section \ref{sec-puzzles-for-smooth}. We prove
a more detailed statement (Theorem \ref{theo-hexp-qft} and
give some consequences in Corollaries
\ref{coro-hexp-max}-\ref{coro-hexp-per}. The first corollary is a
new proof of results in \cite{EE} under an additional assumption.
The second is new.

\subsection{Puzzle and consequences}

At this point, $T$ may be just a continuous self-map of a compact
metric space $M$ together with a finite partition $P$ into subsets
$A$ such that $\bar{A}=\overline{\operatorname{int}(A)}$ and
$T|\bar{A}$ is one-to-one. $\P$ is the set of the interiors of the
elements of $P$.

The puzzle is defined by the refining sequence of ``partitions"
$\P_n$ which are, for each $n\geq1$, the set of \emph{almost
connected components} of the $\P,n$-cylinders, i.e., intersections
of the form $A_0\cap T^{-1}A_1\cap\dots\cap T^{-n+1}A_{n-1}$,
$A_i\in\P_n$. We assume that each $\P_n$ is finite. Their advantage
over the usual connected components is the following key upper bound
on the constraint entropy (to be proved later):

\begin{prop}\label{prop-hC-geo}
For the puzzle $V$ defined by almost connected components of the
$P$-cylinders:
 $$
    h_\C(V) \leq h_\top(T,\partial P) + h_\mult(T,P)
 $$
where $h_\mult(T,P):=\limsup_{n\to\infty} \frac1n\log\mult(P^n)$
with $\mult(Q):=\max_{x\in M} \#\{A\in Q:\overline{A}\ni x\}$.
\end{prop}

We shall show that the puzzle defined in this way by an
entropy-expanding map with a good partition is close to the
original dynamics and also satisfies the remaining assumptions of our
theory.

\begin{remk}\label{rem-puzzle-better}
The above proposition is the counterpart of the upper bound on
minimum left constraint entropy in \cite{QFT}, first claim in the proof of Lemma 7,
p. 385. It is here that we reap the main benefit of the puzzle construction:
we can consider almost connected components of cylinders, instead
of whole cylinders --- thus we get the direct link between the constraint
entropy and the topological entropy of the boundary ``for free", without
having to assume the connectedness of cylinders as in \cite{QFT}, Lemma 7.
\end{remk}

We recall some well-known notions to fix precise definitions and notations.

The \new{coding map} $\gamma_V$ of $(M,T,(\P_n)_{n\geq1})$ (or
just the coding of $V$) is the partially defined map $\gamma:M'\to
X_V$ defined by (i) $M':=\bigcap_{n\geq1}\bigcup_{A\in\P_n}A$;
(ii) $\gamma(x)$ is the unique $y\in X_V$ such that, for all
$n\geq1$, $T^nx\in y_n$. The coding for the usual symbolic
dynamics, simply denoted by $\gamma$, is obtained in this way by
considering the partitions into cylinders of given order:
$\P^1,\P^2,\dots$.

A finite extension of $F:X\to X$ is a skew product over $F$ with finite
fibers, i.e., $G:Y\to Y$ such that $Y\subset X\times\N$, $\#(Y\cap
\{x\}\times\N)<\infty$ for all $x\in X$, and
$G(x,n)=(F(x),\Psi(x,n))$ for some $\Psi:X\times\N\to\N$. We do not
require the cardinality of the fibers to be constant.

A periodic extension of $F:X\to X$ is a map of the form $H:X\times\{0,\dots,p-1\}
\to X\times\{0,\dots,p-1\}$ with, for $0\leq j<p$, $H(x,j)=H(x,j+1)$ and $H(x,p-1)=(F(x),0)$.

\begin{theo}\label{theo-hexp-qft} Let $T:M\to M$ be a $C^\infty$
entropy-expanding map of a compact manifold. Assume that $\P$ is a
good partition and let $(V,i,f)$ be the puzzle obtained by taking
the almost connected components of the $\P,n$-cylinders, $n\geq0$
(see section \ref{sec-dyn-to-puz}). Let  $\gamma_V$ be the coding,

 Then:
 \begin{enumerate}
 \item $\gamma_V$ defines an entropy-conjugacy between $(X_V,F_V)$ and
 $(M,T)$, possibly up to a finite extension: there is
 a Borel finite extension $G$ of $F_V$ and an entropy conjugacy
 of $G$ and $T$ which extends $\gamma_V$;
 \item $h_\C(V)\leq h^{d-1}(T)<h_\top(T)=h_\top(V)$;
 \item $V$ is of quasi-finite type with $h_\wloc(V)=0$;
 \item One can find a determined subpuzzle $V'\subset V$ such that the two previous
 properties still hold and only few periodic orbits are destroyed:
  \begin{equation}\label{eq-det-per}
    \forall N\geq1\; \limsup_{n\to\infty} \frac1n \log \#\{ \xi\in
    i_N(X_V): \xi=F_V^n(\xi) \text{ and } \xi\notin i_N(X_{V'})\} \leq
     h^{d-1}(T).
  \end{equation}
 \end{enumerate}
\end{theo}

Applying Theorems \ref{theo-maxmes} and \ref{theo-class} to $V$
yields a new proof of a slightly weaker version of our result
\cite{EE} about the measures of large entropy of entropy-expanding
maps (we ``lose" here a finite extension):

\begin{coro}\label{coro-hexp-max}
Let $T:M\to M$ be a $C^\infty$ entropy-expanding map. Let $\P$ be
a good partition. Then:
\begin{itemize}
 \item $T$ has finitely many ergodic, invariant probability measure with
maximum entropy;
 \item the natural extension of such maps $T$ are classified up
  to entropy-conjugacy and possibly a period and a finite extension by their
  topological entropy.
\end{itemize}
\end{coro}

Theorem \ref{theo-qft-zeta} applied to $V'$ gives information
about periodic points:

\begin{coro}\label{coro-hexp-per}
In the same setting, let $\eps>0$. Perhaps after replacing $\P$
with a finer good partition, the Artin-Mazur zeta function at level $\P$ of $T$:
 $$
    \zeta_\P(z) := \exp \sum_{n\geq1} \frac{z^n}n
        \#\{\alpha\in \gamma(M'):\sigma^n\alpha=\alpha\}
 $$
is holomorphic on the disk $|z|<e^{-h_\top(T)}$ and has a
meromorphic extension to the larger disk
$|z|<e^{-h^{d-1}(T)-\eps}$. In particular, there exist integers
$p\geq1$ and $m\geq1$ such that for $n\to\infty$ along the
multiples of  $p$:
 $$
    \#\{\alpha\in \gamma(M'):\sigma^n\alpha=\alpha\} \sim
      m e^{n h_\top(f)}
 $$
\end{coro}

\begin{demof}{Corollaries \ref{coro-hexp-max}-\ref{coro-hexp-per}}
Corollary \ref{coro-hexp-max} is a trivial consequence of point 1
of Theorem \ref{theo-hexp-qft} together with Theorems
\ref{theo-maxmes} and \ref{theo-class}.

Corollary \ref{coro-hexp-per} follows similarly from points 3 and 4
of Theorem \ref{theo-hexp-qft} together with Theorem
\ref{theo-qft-zeta} using as the refined finite good partition,
the partition defined by the almost connected components of the
$\P,N$-cylinders where $N=N(\eps)$ is given by Theorem
\ref{theo-qft-zeta}.
\end{demof}

\begin{remk}
(1) If $M$ is one or two-dimensional, then a topological argument
easily shows that each periodic sequence in the coding $
\gamma(M')$ correspond to a periodic point (e.g., using Brouwer
fixed point theorem in connected components of the closure of
cylinders). In higher dimension, one must use the non-uniform
expansion.

(2) The results of Kaloshin \cite{Kaloshin} show that upper bounds
on the number of periodic points can hold for arbitrary maps only
after some identifications.
\end{remk}

In the sequel we prove Theorem \ref{theo-hexp-qft}.

\subsection{Entropy-conjugacies}

\begin{lem}\label{lem-h-conj}
Let $T:M\to M$ be an entropy-expanding map with a good partition
$\P$. Then the puzzle defined by the almost connected components of
cylinders has the same entropy as $T$. More precisely, the coding
$\gamma_V$ defines an entropy-conjugacy between $T$ and a Borel
finite extension of $F_V$.
\end{lem}

To prove this, we use a common extension $X_V\ltimes M$ of the
puzzle and of $T$ defined as:
 $$
    X_V\ltimes M = \overline{ \{(v,x)\in X_V\times M:\forall
    n\geq0\; x\in v_n\} }
 $$
endowed with the map $F_V\ltimes T$ which is just the restriction
of the direct product. Let $\pi_1$, resp. $\pi_2$, be the
projection $X_V\ltimes M\to X_V$, resp. $X_V\ltimes M\to M$.

\medbreak\noindent $\bullet$ We claim that $F_V\ltimes T$ and $T$
are entropy-conjugate. Observe that, the partition being good for
$T$, no point returns infinitely many times to $\partial P$. Hence
$\partial P$ has zero measure w.r.t. any $T$-invariant probability
measure. The same is true for $\pi_2^{-1}(\partial P)$. Hence
$(v,x)\mapsto x$ is an isomorphism w.r.t. any invariant
probability measure, proving the claim. In particular,
$h_\top(F_V\ltimes T)=h_\top(T)$ by the variational principle.

\medbreak\noindent $\bullet$ We claim that $F_V\ltimes T$ and
$F_V$ are entropy-conjugate, perhaps after replacing the latter
$F_V$ by a Borel finite extension. As the extension is continuous
and compact, any invariant probability measure of $F_V$ can be
lifted to $F_V\ltimes T$. We have to show that, given a large
entropy measure of $F_V$ (1) there are only finitely many ergodic
lifts $\hat\mu$; (2) for each such $\hat\mu$, $\pi_1: (X_V\ltimes
M,\hat\mu) \to(X_V,\pi_1\hat\mu)$ is a finite extension.

We first prove point (2). We can assume $\hat\mu$ to be an arbitrary
$F_V\ltimes T$-invariant and ergodic probability measure with
$h(F_V\ltimes T,\hat\mu)= h(T,\pi_2\hat\mu)>h^{d-1}(T)$. Let $\mu=(\pi_2)_*\hat\mu$ and
$\nu=(\pi_1)_*\hat\mu$. $\mu$ is a $T$-invariant ergodic measure
satisfying $h(T,\mu)>h^{d-1}(T)$. By \cite{EE} this implies that
$\mu$ has only strictly positive Lyapunov exponents, hence, by
\cite{coding}, $\pi_1:(X_V\ltimes M,\hat\mu)\to(X_V,\nu)$ is a
finite extension. This proves point (2).

We prove point (1) following \cite{coding}. Assume by
contradiction that there exists infinitely many distinct ergodic
lifts $\hat\mu_1,\hat\mu_2,\dots$ of some ergodic and invariant
probability measure $\mu$ of $F_V$. We can assume that $\hat\mu_n$
converges to some $\hat\mu_*$. As $\pi_2$ is continuous,
$\hat\mu_*$ is also a lift of $\mu$ and so are almost all of its
ergodic components.  They project on $M$ to ergodic invariant
probability measures with positive Lyapunov exponents. As
explained in \cite{coding}, this implies that for each such
ergodic component $\hat\nu$, for $\hat\nu$-a.e. $(v,x)$, there
exists a ball $B$ around $x$ in the fiber which contains no
generic point wrt any measure distinct from $\hat\nu$. It follows
that there are only countably (or finitely) many ergodic
components. Thus, there exists an ergodic component of
$\hat\mu_*$, such that the union of these fibered neighborhood has
positive $\hat\mu_*$-measure. Hence it has positive measure for
$\hat\mu_n$ for $n$ large. But this implies that
$\hat\mu_n=\hat\mu_*$, a contradiction. Point (1) is proven and
the claim follows.

\medbreak

The above two claims prove the lemma.

\subsection{Constraint entropy}

Before proving Proposition \ref{prop-hC-geo} which will imply claim 2
of Theorem \ref{theo-hexp-qft}, we give a geometric necessary
condition for the irreducibility of puzzle pieces.

\begin{lem}\label{lem-geo-irr}
Let $(V,i,f)$ be a puzzle generated by the almost connected
components of the cylinders of a partition $\P$.

Let $v\in V$ and let $A$ be the unique the element of $\P$
containing $v$,
 $$
    \overline{f(v)}\cap\partial T(A)=\emptyset \implies v \text{
    is $f$-reducible}.
 $$
\end{lem}

\begin{demo}
Assume $\overline{f(v)}\cap\partial T(A)=\emptyset$. $v$ is an
almost connected component of $A\cap T^{-1}(f(v))=
(T|\bar{A})^{-1}(f(v))$. By the assumption, this last set is
uniformly homeomorphic to $f(v)$, hence is almost connected.
Therefore it is equal to $v$.

This shows that $v$ is uniquely determined by $f(v)$ and
$A=i_1(v)$ ({\it a fortiori} $i(v)$), verifying condition (2) of
reducibility.

Consider now (*) $f:\itree(v)\to\itree(f(v))$. Observe that for
any $w\in\itree(f(v))$, $w\subset f(v)$. Hence, $\overline{w}\cap
\partial T(A)=\emptyset$. The reasoning for the uniqueness of
$v$ shows that the map (*) is one-to-one: $f(u)=f(u')$ implies
that $T(u)$ and $T(u')$ are both almost connected subsets of
$f(u)=f(u')$, so they must be equal.

For $w\in\itree(f(v))$,  $u=(T|A)^{-1}(w)\in\itree(v)$ satisfies
$w=f(u)$. Hence the map (*) is onto and therefore an isomorphism,
proving condition (1) of reducibility.
\end{demo}

\begin{demof}{Proposition \ref{prop-hC-geo}}
Let $r>0$ and $\eps>0$. Let $\Sigma_n$ be an arbitrary
$(r,n)$-separated subset of $\C_n$, the set of irreducible pieces
of order $n$. Recall that there exists some $L=L(r)$, such that,
for all $n\geq L$, $x,x'\in V_n$ are $(r,n)$-separated then there exists
some $0\leq k<n-L$ such that $(f^kx)_L\ne (f^kx')_L$.

We are going to bound the cardinality of $\Sigma_n$ by $e^{(h_\top(T,\partial\P)+h_\mult(T,P)+2\eps)n}$.

Let $\rho>0$ be smaller than the distance between any two almost
connected component of any $L$-cylinder (there are only finitely
many of them, $L$ being fixed, and the distance between any two of
them is positive as we are considering {\bf almost} connected
components). For all integers $n$ large
enough, $\mult(\P^n)< e^{(h_\mult(T,P)+\eps) n}$ and
$r(\rho/2,n,\partial\P)\leq e^{(h_\top(T,\partial\P)+\eps)n}$.

Let $S_n$ be a minimum $(\rho/2,n)$-spanning subset
of $\partial P$. To every $v\in \Sigma_n$, associate a point
$x=x(v)\in S_n$ such that $d(T^kv,T^kx)<\rho/2$ for all $0\leq
k<n$ ($T^kv$ is a subset of $M$). This is possible since $\bar v\cap\partial P\ne\emptyset$ by
Lemma \ref{lem-geo-irr}.

The map $x:\Sigma_n\to S_n$ is at most $(\#V_L)^L\cdot e^{(h_\mult(T,P)+\eps)
n}$-to-1. Indeed, assume that there exists $x\in S_n$ with more
than this number of pre-images. We can find a set of $e^{(h_\mult(T,P)+\eps)
n}$ pre-images, all with the same $((f^kv')_L)_{n-L\leq k<n}$. As
$e^{(h_\mult(T,P)+\eps) n}> \mult(\P^n)$, two of these,
say $v$ and $v'$, must almost connected components of the same $n,\P$-cylinder. But
then $d(T^kv,T^kv')\leq d(T^kv,T^kx)+d(T^kx,T^kv')<\rho$ implies
$(f^kv')_L=(f^kv)_L$ for all $0\leq k<n-L$, contradicting the separation
assumption. Therefore $\#\Sigma_n\leq (\#V_L)^L\cdot e^{(h_\mult(T,P)+\eps) n} \# S_n$
and
 $$
    h_\C(F_V) \leq h_\top(T,\partial\P)+h_\mult(T,P)+2\eps,
 $$
with arbitrary $\eps>0$, proving the claim.
\end{demof}

\subsection{Determinacy}

We turn to determinacy. The delicate point here is that it is possible (though exceptional)
that $u\fred v$ in the absence of the geometric property of Lemma
\ref{lem-geo-irr}, because of the following phenomenon.

A puzzle piece $v\in V$ is \new{trivial} if there exists $k\geq1$
such that for every $w\in\itree(v)$,
 \begin{equation}\label{eq-bad-w}
   w\cap T(\partial\P^k)\ne\emptyset.
 \end{equation} The trivial subset of $V$ is
the smallest subset $V^0$ of the puzzle such that:
 \begin{itemize}
  \item $V^0$ contains all trivial pieces;
  \item if $f(v)\in V^0$ then $v\in V^0$.
 \end{itemize}
Observe that $V\setminus V^0$ equipped with the restrictions of $i$
and $f$ is a puzzle as $f(V\setminus V^0)\subset V\setminus V^0$ (by
definition) and $i(V\setminus V^0)\subset V\setminus V^0$ (as $\itree(i(v))\supset
\itree(v)$).

\begin{lem}
Let $(V,i,f)$ be the puzzle defined by a dynamical system $T:M\to
M$ as in Proposition \ref{prop-hC-geo}.
The non trivial puzzle $V':=V\setminus V^0$ is {determined}.

Assume additionally that $h_\C(V)<h_\top(V)$. Then the obvious
injection $i:X_{V'}\to X_V$ is an entropy-conjugacy and the
approximate periodic points of the two systems satisfy the
estimate (\ref{eq-det-per}) of Theorem \ref{theo-hexp-qft}:
 $$
   \forall N\geq 1\; \limsup_{n\to\infty} \frac1n\log
     \#\{\xi\in i_N(X_V)\setminus i_N(X_{V'}) : \sigma^n\xi=\xi \}
        \leq h_\C(V).
 $$
\end{lem}

\begin{demo}
To prove the determinacy, we consider $v,v',w\in V'$ such that
$i_1(v)=i_1(v')=:A\in\P$ and $v,v'\fred^1 w$. Assume by
contradiction that $v\ne v'$. $T|\bar{A}$ is a homeomorphism hence
$v,v'$ are disjoint almost connected components of $(T|\bar{A})^{-1}(TA\cap
w)$. $v,v'\fred w$ implies that:
 $$
   \{f(u):u\in\itree(v)\} = \{f(u'):u'\in\itree(v')\} =\itree(w).
 $$
Hence, every $t\in\itree(w)$ is an almost connected set containing
both disjoint sets $T(v)$ and $T(v')$. Therefore $t\cap \partial T(v)
\ne\emptyset$, so that $t\cap T(\partial\P^{|v|})\ne\emptyset$. Thus $w$ is
trivial, the sought-for contradiction.

\medbreak

We now let $\mu$ be an ergodic $F_V$-invariant probability measure
such that, for some $v\in V^0$, $\mu([v]_V)>0$. By invariance of $\mu$,
$\mu([w]_V)>0$ for a trivial $w=f^n(v)$ with $0\leq n<|v|$.
Now, $x\in[w]_V$ implies that $x_{|w|}=w$ and $x_n\cap T(\partial\P^{\ell})\ne\emptyset$
for all $n\geq|w|$ and some fixed, minimal $\ell\geq0$. Therefore $t:=F_V^\ell(x)$ satisfies
$t_m\cap T(\partial P)$ for all $m\geq |w|-\ell$. The reasoning in
the proof of Proposition \ref{prop-hC-geo} implies:
 $$
    h(F_V,\mu)\leq h_\top(F_V,[w]_V) \leq h_\C(V) < h_\top(V),
 $$
proving the entropy-conjugacy.

Consider now some periodic sequence $\xi\in i_N(X_V)\setminus
i_N(X_{V'})$. Hence $\xi=i_N(x)$
with $x_m\in V^0$ for some $m\geq0$. Thus $f^n(x_m)$ is trivial.
We may assume $n=0$ by shifting to another point of the same periodic orbit).
Therefore $x_p\cap T(\partial\P^k)\ne\emptyset$ for some $k$
and all $p\geq m$. As above, it follows that $f^{j}(x_{p+j})\cap T\partial\P
\ne\emptyset$ for all $p\geq m$ and some $j$ which can be assumed
to be fixed and then $0$.
The claimed bound on the number of periodic points follows.
\end{demo}

\subsection{W-local Entropy}

We prove the third point of Theorem \ref{theo-hexp-qft}:

\begin{lem}
If $V$ is a \sQFT puzzle which is determined then, for all
invariant and ergodic probability measures $\mu$ on $X_V$ with
$h(F_V,\mu)>h_\C(V)$,
 $$
    h(F_V,\mu)=h(F_V,\mu,\eps_*)
 $$
So in particular, $h_\wloc(V)=0$.
\end{lem}

\begin{demo}
Any ergodic invariant probability measure on $X_V$ with entropy $>h_\C(V)$
can be lifted to an isomorphic $\hat\mu$ on $\Sigma(\D)$ by
Theorem \ref{theo-struct}. $h(\sigma^{-1},\hat\mu)=
h(\sigma,\hat\mu) =h(F_V,\mu)$ can be bounded by the growth rate
of the number paths on $\D$ ending at any fixed vertex $v_*\in\D$
with $\hat\mu([v]_{\Sigma(\D)})>0$. But those paths are uniquely
determined by their $i_1$-projection as $V$ is determined (Lemma
\ref{lem-determined}). Thus,
$h(F_V,\mu)=h(\sigma,i_1(\mu))=h(F_V,\mu,1/2)$.
\end{demo}

\begin{coro}
Let $V$ be a \sQFT puzzle with a subpuzzle $V'$ which is
determined. Assume that the inclusion $X_{V'}\to X_v$ is a
conjugacy with respect to all ergodic invariant
probability measures with entropy $>h_\C(V)$. Then $h_\wloc(V)=0$.
\end{coro}

\begin{demo}
By the previous lemma, $h_\wloc(V')=0$. Let us see that this
property carries over to $V$.

Let $\mu$ be an ergodic invariant probability measure of $X_V$
with $h(F_V,\mu)>h_\C(V)$. Hence, it can be identified to an
invariant measure $\mu'$ of $F_{V'}$. Therefore $h(F_{V'},\mu')=
h(\sigma,i_N(\mu'))$ for some integer $N\geq1$. But $V'\subset V$
hence one can define almost everywhere $i_N:X_V\to i_N(X_{V'})$
and check that $i_N(\mu)$ and $i_N(\mu')$ are isomorphic so that
$h(F_{V},\mu)= h(\sigma,i_N(\mu))$, proving the claim.
\end{demo}

\appendix

\section{Varying Radius of Meromorphy}

\begin{defi}
Denote by $\MR(f)$ the radius of meromorphy of a formal power
series $f$. It is zero if the radius of convergence of $f$,
$\rho(f)$, is zero. Otherwise it is the supremum of the radiuses
$r$ of the disks $D(r)$ centered at zero for which there exists a
rational function $F(z)$ such that $f(z)/F(z)$ can be extended to
a holomorphic and non-zero function on $D(r)$.
\end{defi}

\begin{fact}
There exists a countable oriented SPR graph $G\ni a,b$ such that
$\MR(\zeta^G_a)\ne \MR(\zeta^G_b)$.
\end{fact}

We found this example after an illuminating discussion with O.
Sarig.

\medbreak

Before giving our construction, we recall some basic tools. The
main tool here is the notion of a loop graph (or petal graph in
the terminology of B. Gurevi\v{c}). These graphs have a
distinguished vertex and an arbitrary number of first return
loops\footnote{That is, sequences $v_0\to^{e_1}
v_1\to\dots\to^{e_n} v_n$ where the vertices $v_i$'s and edges
$e_i$'s are distinct except for $v_0=v_n$ which is the
distinguished vertex.} of each length based at the distinguished vertex,
but distinct first return loops are disjoint except for the
distinguished vertex. Such graphs are completely described by
their \new{first return series} $f(z):=\sum_{n\geq1} f_nz^n$ where
$f_n$ is the number of first return loops of length $n$ (based at the
distinguished vertex). It is well-known that the local zeta
function at the distinguished vertex is
 $$
   \zeta^G_*(z) = \frac{1}{1-f(z)} = \sum_{n\geq1} \ell_n z^n
 $$
where $\ell_n$ is the number of loops of length $n$ based at the
distinguished vertex\footnote{These loops may go several times
through the distinguished vertex.}.

\medbreak

We now give the construction. We consider two disjoint loop graphs
defined by first return series $a(z):=\sum_{n\geq1} a_nz^n$ and
$b(z):=\sum_{n\geq1} b_nz^n$. We call their respective
distinguished vertices $a$ and $b$.

We define a new graph $G$ by taking the disjoint union of:
 \begin{itemize}
  \item the two preceding loop graphs;
  \item a set of disjoint paths from $a$ to $b$ described by
a series $s(z):=\sum_{n\geq1} s_nz^n$ (there are $s_n$ simple
paths of length $n$ from $a$ to $b$ and these are disjoint);
  \item a set of simple
paths (i.e., injective as maps) from $b$ to $a$ described by a series $t(z):=\sum_{n\geq1}
t_nz^n$.
 \end{itemize}

\begin{claim}
The first return series of $G$ at $a$ is:
 $$
   \hat a(z) = a(z) + \frac{ s(z)t(z) }{1- b(z)}.
 $$
\end{claim}

Indeed, any first return loop at $a$ in $G$ is exactly in
one of the following classes:
 \begin{itemize}
  \item the first return loops in the loop graph $a$;
  \item the concatenations of a transition from $a$ to $b$, a (not necessarily first
  return) loop at $b$, a transition from $b$ to $a$.
 \end{itemize}

Fix $b(z)=2z^2$ (so the associated Markov shift is the set of all
infinite concatenations of the two words of length $2$, say $b0$
and $b1$).

Let $\tau(z):=\sum_{n\geq1} \tau_n z^n=s(z)t(z)$. We arrange
it so:
 \begin{itemize}
   \item $\tau_0=\tau_1=0$, $\tau_n=0$ or $1$;
   \item $|z|=1$ is the natural boundary of $\tau$.
 \end{itemize}
By the P\'olya-Carlson theorem \cite{Complex}, the last condition
is equivalent to $\tau_n$ not being eventually periodic.
This can be obtained by taking $\{n:s_n=1\}$ and $\{n:t_n=1\}$ to
be disjoint subsets of $2,2^2,2^3,\dots$ satisfying the aperiodic
condition above. It follows that $q(z):=\tau(z)/(1-b(z))$ satisfies:
 \begin{enumerate}
  \item $q_0=0$, $0\leq q_n\leq 2^{n/2+1}\leq 5^n$
  \item $|z|=1$ is the natural boundary of $q$.
 \end{enumerate}
Now, set $a_0=0$ and, for $n\geq1$:
 $
   a_n := 5^n - q_n \geq 0.
 $
We have: $\hat a_n=a_n+q_n=5^n$ for $n\geq 1$, $\hat a_0=0$. Hence
$\hat a(z)=5z/(1-5z)$ and
 $$
   \zeta^G_a(z) = \frac1{1-\hat a(z)} = \frac{1-5z}{1-10z}
 $$
is a rational function. In particular, $\MR(\zeta^G_a) =\infty$.
On the other hand,
 $$
   \hat b(z) = b(z) + \frac{\tau(z)}{1-a(z)}
    = 2z^2 + \frac{(1-2z^2)q(z)}{1-\frac{5z}{1-5z}+q(z)}
    = 2z^2 + \frac{(1-2z^2)(1-5z)}{1+(1-10z)/q(z)}
 $$
Therefore $\zeta^G_b(z) = 1/(1-\hat b(z))$ has meromorphy radius:
 $
    \MR(\zeta^G_b) = \MR(\hat b) = \MR(q) = 1
 $
and $$\MR(\zeta^G_b) < \MR(\zeta^G_a)$$ as claimed.

Observe that $h(G)=\log 10$ and $h_\infty(G)=\log 5$. Hence
$h_\infty(G)<h(G)$ and $G$ is SPR as claimed, finishing the
construction.

\section{Good Partitions for Almost All Couplings}\label{appendix-partition}

We consider the following, convenient family of coupled maps.
For $(a,b,c)\in\R^3$, we let
 $$
    F_{a,b,c}(x,y)=\left(a(1-4x^2)+cy^2-\frac12,b(1-4y^2)+cx^2-\frac12\right).
 $$
For $(a,b,c)\in\Omega=\{(a,b,c)\in(0,1)^3: c<4-4\max(a,b) \}$,
$F_{a,b,c}(Q)\subset Q$ for $Q:=[-\frac12,\frac12]^2$.
There is a natural partition $\P$ into four elements:
 $$
  Q_{\eps_1,\eps_2}:=\{(x,y)\in Q:\eps_1x>0,\; \eps_2y>0\}
    \qquad (\eps_1,\eps_2)\in\{-1,1\}^2
 $$
according to the signs of $x$ and $y$. Most of the properties of a good partition are obvious for this convenient family:

Indeed, $F_{a,b,c}|\bar Q_{\eps_1,\eps_2}$ is obviously one-to-one.  The boundary of the partition is:
$$
  \partial\P=[-1/2,1/2]\times\{0\}\cup\{0\}\times [-1/2,1/2]
   \cup\partial Q.
$$
$\partial \P$ is obviously the image of a compact subset of $\R$ by a $C^\infty$ map. The semi-algebraic nature of both $F_{a,b,c}$ and $\partial\P$ implies that each cylinder has indeed finitely many connected components and therefore finitely many almost connected components.

To conclude, we show that, after discarding countably many hypersurfaces in the parameter space, there is a constant such that:
 \begin{equation}\label{eq:three}
   \forall(x,y)\in Q\;  \#\{k\geq0: F_{a,b,c}^k(x,y)\in\partial\P\} \leq 2.
 \end{equation}
We prove that for each $0<n<m$, there exists a hypersurface containing all the parameters $(a,b,c)$ such $F_{a,b,c}^n(x,0)$ and $F_{a,b,c}^m(x,0)$ are both in $\{0\}\times[-1/2,1/2]$ for some  $x\in[-1/2,1/2]$. The cases involving other pieces of $\partial\P$ are similar
and together they imply eq. (\ref{eq:three}).

Observe that $F_{a,b,c}^k(x,0)\in\{0\}\times[-1/2,1/2]$ is equivalent to
 $$
   P_{k,a,b,c}(x)=0
 $$
for some polynomials in $x$ whose coefficients are themselves polynomials in $a,b,c$.

The degrees of $P_{n,a,b,c}$ and $P_{m,a,b,c}$ are fixed, say $p$ and $q$, outside of an algebraic hypersurface . Hence the parameters we are considering are such that the $(p+q)\times(p+q)$ resultant of the two polynomials $P_{n,a,b,c}$ and $P_{m,a,b,c}$ is zero: these parameters satisfy a polynomial equation. This equation is not trivial as it is not satisfied for $a=b=1$, $c=0$.
Indeed, $F_{1,1,0}(x,y)=(\frac12-4x^2,\frac12-4y^2)$ so the $x$-coordinate can take the value $0$ only once in an orbit (the subsequent values are then $1/2,-1/2,-1/2,\dots$).

\end{document}